\documentclass{amsart}
\usepackage{graphics}
\usepackage{maple2e}
\newlength{\itemwidth}
\newenvironment{imini}{\addtolength\itemwidth{-\leftmargin}\addtolength\itemwidth{-\labelwidth}\\\begin{minipage}{\itemwidth}}{\end{minipage}\\}

\numberwithin{equation}{section}
\newtheorem{Lem}{Lemma}[section]
\newtheorem{Prop}[Lem]{Proposition}
\newtheorem{Cor}[Lem]{Corollary}

\theoremstyle{definition}
\newtheorem{Def}[Lem]{Definition}
\theoremstyle{remark}
\newtheorem{Rem}[Lem]{Remark}
\newtheorem{Expl}[Lem]{Example}
\newcommand\Hit{\rhd}
\newcommand\hit{\rightharpoonup}

\newcommand\hitby{\leftharpoonup}
\newcommand\Hitby{\lhd}
\renewcommand\o{\otimes}
\DeclareMathOperator\id{\operatorname{id}}
\newcommand\op{{\operatorname{op}}}
\newcommand\cop{{\operatorname{cop}}}

\newcommand\sw[1]{{}_{(#1)}}
\newcommand\pow[1]{{}^{[#1]}}
\newcommand\ol{\overline}
\newcommand\inv{^{-1}}
\renewcommand\epsilon\varepsilon

\newcommand\C{\mathbb C}
\newcommand\Z{\mathbb Z}
\newcommand\N{\mathbb N}
\newcommand\QQ{\mathbb Q}
\newcommand\bas[1]{p_{{#1}}}
\newcommand\neut{1}
\newcommand\TPS{T}
\newcommand\TPM{T'}
\newcommand\tpd{t}
\newcommand\tpn{t'}
\newcommand\Ker{\operatorname{Ker}}
\newcommand\rank{\operatorname{rank}}
\newcommand\ch{\operatorname{char}}
\makeatletter
\def\namelabel#1#2{\@bsphack
  \protected@write\@auxout{}%
         {\string\newlabel{#1.nme}{{#2}{#2}}}%
  \@esphack}
\def\nmlabel#1#2{\label{#2}\namelabel{#2}{#1}}
\newcommand\nmref[1]{\ref{#1.nme}\ \ref{#1}}
\makeatother

\begin{document}
\title{Hopf powers and orders for some bismash products}
\author{Rachel Landers}
\author{Susan Montgomery}
\address{Department of Mathematics, USC, Los Angeles, CA 90089-1113,
USA}
 \email{smontgom@math.usc.edu, landers@usc.edu}
\author{Peter Schauenburg}
\address{Department of Mathematics, USC, Los Angeles, CA 90089-1113, USA,
and Mathematisches Institut der Universit\"at M\"unchen,
Theresienstr.~39, 80333~M\"unchen, Germany}
\email{schauen@mathematik.uni-muenchen.de} \subjclass{16W30}
\thanks{}
 \keywords{}

\maketitle
\section*{Introduction}
Let $H$ be a Hopf algebra over the field $k$. We study the $n$-th
Hopf power map of $H$, a linear endomorphism $[n]\colon H\to H$.
If $H=kG$ is a group algebra, then $[n]$ is the linear extension
of the $n$-th power map on $G$. If $H$ is commutative, and thus
represents an affine group scheme $G$, the $n$-th Hopf power map
represents the endomorphism of the affine scheme $G$ given by
taking the $n$-th power of all elements of the group $G(R)$ for
each $k$-algebra $R$. In this sense the Hopf power map for
commutative Hopf algebras is implicit in work of Gabriel
\cite{Gab:EISGGF}; explicitly it seems to appear first in a paper
by Tate and Oort \cite{TatOor:GSPO}. Kashina
\cite{Kas:OAHAHHYD,Kas:GPMHA} proposed studying the same map for
general (non-commutative and non-cocommutative) Hopf algebras. If
$H$ is commutative and represents the affine group scheme $G$,
then the $n$-th Hopf power map for $H$ is trivial if and only if
the $n$-th power map for each of the discrete groups $G(R)$ is
trivial. Thus the smallest number such that $[n]$ is trivial can
be viewed as a sort of joint exponent for all the groups $G(R)$,
or the exponent of the algebraic group $G$. Kashina reproves a
result of Gabriel \cite[Prop. 8.5]{Gab:EISGGF} using more
Hopf-algebraic techniques: For an $n$-dimensional commutative Hopf
algebra the $n$-th Hopf power map is trivial, i.e. the exponent of
a finite group scheme over a field divides its order. On a more
elementary level, the exponent of a finite group $G$ is also the
least positive integer for which the $n$-th Hopf power map on the
group algebra $kG$ is trivial.

For a general finite-dimensional Hopf algebra $H$, Kashina asks
whether $[n]$ is trivial for some finite number $n$, preferably a
divisor of $\dim H$. She proves several positive results, most
notably that $[n]$ is trivial for the Drinfeld double of the group
algebra of a group of exponent $n$; more generally, $[n]$ is
trivial for $D(H)$ if it is for $H$ and $H^\op$.

Etingof and Gelaki \cite{EtiGel:EFDHA} mint the obvious term
exponent of $H$ for the least $n$ such that $[n]$ is trivial. They
point out that Kashina's definition of $[n]$ is only suitable for
involutive Hopf algebras, and define the right version for the
general case. Among the various important results of
\cite{EtiGel:EFDHA} is that the exponent of a semisimple Hopf
algebra $H$ over a field of characteristic $0$ is indeed finite
and divides $(\dim H)^3$, and the exponent is invariant under
Drinfeld twists.

In the present paper we study the behavior of individual elements
of $H$ under the Hopf power maps. For which $n$ are there
nontrivial elements of $H$ whose $n$-th Hopf power is trivial?
What are the possible Hopf orders of elements of $H$ (where the
Hopf order of $h$ is the least positive integer $n$ such that the
$n$-th power of $h$ is trivial)? Since the $n$-th Hopf power map
for a group algebra $kG$ is the linearization of the $n$-th power
map of $G$, these questions generalize (or linearize) questions
about groups whose answers are well-known. Although we will see
that the answers for group algebras deviate somewhat from the
answers for groups, they are still quite easy and reasonable.

Before turning to more interesting examples, we propose a
refinement of the above questions. The space $\TPS_n(H)$ of all
elements of $H$ whose $n$-th Hopf power is trivial is a linear
subspace of $H$. The dimensions $\tpd_n(H)$ of all these spaces
$\TPS_n(H)$, as well as the dimensions $\tpd_{m,n}(H)$ of their
pairwise intersections, are isomorphism invariants of $H$; in
particular they tell us for which $n$ there exist nontrivial
elements whose $n$-th power is trivial, and for which $n$ there
exist elements of Hopf order $n$. If $H$ is a Kac Hopf algebra,
the numbers $\tpd_{m,n}$ have certain symmetries (apart from the
obvious $\tpd_{m,n}=\tpd_{n,m}$).

We study the power maps, possible Hopf orders, and the dimensions
$\tpd_{m,n}$ for the (semisimple) bismash product Hopf algebras
obtained from a factorizable group. Since explicit calculations
for these examples are quite involved, we rely heavily on computer
help; for most of our examples we use Maple to first compute the
matrices of the Hopf power maps with respect to a suitable basis,
and then to do computations with these matrices.

Here is a brief summary of some of our findings: A Hopf algebra
$H$ can have elements whose Hopf order is a prime that does not
divide the dimension of $H$. The numbers $\tpd_{m,n}$ are not
invariant under Drinfeld twists. While it is easy to check that
the numbers $\tpd_n$ are the same for $H$ and its dual $H^*$, the
numbers $\tpd_{m,n}(H)$ and $\tpd_{m,n}(H^*)$ can be different.
Also, the Hopf orders of elements that occur in $H$ can be
different from the Hopf orders of elements that occur in $H^*$.
(We should note at this point that the possible Hopf orders of
elements of a group algebra and its dual are the same; computer
experiments show that $\tpd_{m,n}(\QQ G)=\tpd_{m,n}(\QQ^G)$ for
all of the few groups that we have checked.)
\subsection*{Acknowledgements}
Work on this paper was started while the third author was visiting
the University of Southern California --- in particular the second
author --- during the summer of 2003. He thanks the Department of
Mathematics at USC for a warm welcome and excellent working
conditions, and the \emph{Deutsche Forschungsgemeinschaft} for
generous support by a Heisenberg Fellowship. The second author was
supported by NSF grant DMS-0100461. The first author was supported
by a \emph{Women in Science and Engineering} (WiSE) undergraduate
research grant at USC.
\section{Definitions, and the (co)commutative case}
Throughout the paper we work over a fixed base field $k$. Hopf
algebras, vector spaces, and tensor products are understood to be
over $k$. We use Sweedler notation in the form $\Delta(h)=h\sw 1\o
h\sw 2$ for comultiplication.
\begin{Def}
Let $H$ be a Hopf algebra. For $n\in\Z$, the \emph{$n$-th Hopf
power map}  $[n]=[n]_H\colon H\to H$ is the $n$-th convolution
power of the identity. Thus, writing $[n](h)=:h\pow n$, we have
$h\pow n=h\sw 1h\sw 2\cdot\dots\cdot h\sw n$ for $n>0$, $h\pow
0=\epsilon(h)1$, and $h\pow n=S(h\sw 1)S(h\sw 2)\cdot\dots\cdot
S(h\sw {-n})$ for $n<0$. The element $h\pow n$ is called the
\emph{$n$-th Hopf power} of $h$.

We will say that the $n$-th Hopf power of $h$ is \emph{trivial} if
$h\pow n=\epsilon(h)1$. The \emph{Hopf order} of $h$ is the least
positive integer $n$ such that the $n$-th Hopf power of $h$ is
trivial. The \emph{exponent} of $H$ is the least positive integer
$n$ such that the $n$-th Hopf power of every element of $H$ is
trivial.
\end{Def}
Etingof and Gelaki \cite{EtiGel:EFDHA} have shown that the above
definition of Hopf powers and the exponent is suitable only for
involutive Hopf algebras. For the general case they give a
better-behaved definition involving the square of the antipode.
The most interesting class of Hopf algebras for the investigation
of Hopf powers and the exponent still seems to be that of
semisimple Hopf algebras over a field of characteristic zero. We
will therefore stick to the above definition, even though we will
not always make those assumptions on the Hopf algebras in
consideration.

 Most
results will require the Hopf algebra $H$ to have finite exponent.
\begin{Rem}\nmlabel{Remark}{genpowrem}
For all $m,n\in\Z$ we have $h\pow{m+n}=h\sw 1\pow mh\sw 2\pow n$,
and $S(h\pow n)=S(h)\pow n$.

In particular, $h\pow{n+ek}=h\pow n$ for all $n,k\in\Z$, if $e$
denotes the exponent of $H$.
\end{Rem}
\begin{Def}
  Let $H$ be a Hopf algebra. The \emph{$n$-th trivial power space} of $H$ is
  $$\TPS_n:=\TPS_n(H):=\{h\in H|h\pow n=\epsilon(h)1\}=\Ker([n]-\eta\epsilon).$$
The \emph{$n$-th trivial power dimension} of $H$ is
$\tpd_n:=\tpd_n(H):=\dim\TPS_n(H)$. We put
$\TPS_{m,n}:=\TPS_{m,n}(H):=\TPS_m\cap\TPS_n$, and
$\tpd_{m,n}:=\tpd_{m,n}(H):=\dim(\TPS_{m,n})$.
\end{Def}
Since $1\pow n=1=\epsilon(1)\cdot 1$ for all $n$, all spaces
$\TPS_{m,n}$ and all dimensions $\tpd_{m,n}$ are nonzero.
Nontrivial elements with trivial $n$-th Hopf power exist whenever
$\tpd_n>1$. Over an infinite field, the numbers $\tpd_{m,n}$ also
determine whether elements of a given Hopf order exist:
\begin{Prop}\nmlabel{Proposition}{ordersbydim}
  Assume that $H$ is a finite dimensional Hopf algebra over an
  infinite base field $k$.
  Then the following are equivalent
  \begin{enumerate}
    \item $H$ contains an element of Hopf order $n$.
    \item $\TPS_n(H)\not\subset\TPS_m(H)$ for all $m<n$.
    \item $\tpd_n(H)>\tpd_{m,n}(H)$ for all $m<n$.
  \end{enumerate}
\end{Prop}
\begin{proof}
  By definition there is an element of Hopf order $n$ if and only if
  $\TPS_n(H)$ contains an element not in $\bigcup_{m<n}\TPS_m(H)$,
  or if $\TPS_n(H)$ properly contains
  $\bigcup_{m<n}(\TPS_m(H)\cap\TPS_n(H))$.
  But a vector space over $k$ cannot be written as the union of
  finitely many proper subspaces (see for example \cite[Sec.~4.2, Ex.~21]{Her:TA}),
  so this is equivalent to
  the condition that $\TPS_n(H)$ properly contains each of the
  spaces $\TPS_m(H)\cap\TPS_n(H)$ for $m<n$. This is equivalent on
  one hand to the proper inequality of dimensions, and on the
  other hand to $\TPS_n(H)$ not being contained in $\TPS_m(H)$.
\end{proof}
\begin{Cor}
  Let $H$ be a Hopf algebra with finite exponent $e$ over an infinite
  base field $k$. Then $H$ contains an element of Hopf order $e$.
\end{Cor}
\begin{proof}
  By definition, $e$ is the least integer such that $H=\TPS_e(H)$.
  If there is no element of Hopf order $e$, the preceding result shows
  that there is $m<e$ with $H=\TPS_e(H)\subset\TPS_m(H)\subset H$,
  a contradiction.
\end{proof}

If $H$ is semisimple over the complex numbers, then Kashina,
Sommerh\"auser, and Zhu in fact show that an integral in $H$ has
Hopf order $e$, see the end of Section 3 in \cite{KasSomZhu:HFSI}.
This is easy to check explicitly for a group algebra $kG$, and
also for its dual (cf. \nmref{kGex} and \nmref{kG*ex}
\eqref{kG*int})

In the rest of this section we collect the facts on Hopf powers
and the trivial power spaces for commutative or cocommutative Hopf
algebras and compare to the more familiar notions for groups.
\begin{Def}
  Let $H$ be a Hopf algebra. We will say that $H$ satisfies the \emph{power rule} if
\begin{equation}\label{powrule}
  (h\pow m)\pow n=h\pow {mn}
\end{equation}
holds for all $h\in H$, $m,n\in\Z$.
\end{Def}
\begin{Lem}
  Any cocommutative or commutative Hopf algebra satisfies the power rule.
\end{Lem}
This observation was made by Tate and Oort \cite{TatOor:GSPO} for
commutative Hopf algebras. A proof for positive exponents avoiding
the language of group schemes is in Kashina's paper
\cite{Kas:GPMHA}. Negative exponents are also easy to cover, more
so perhaps after \nmref{powbem} below.

Clearly the power rule also holds for the tensor product of two
Hopf algebras that satisfy the power rule. It also carries over to
dual Hopf algebras, quotient Hopf algebras, and Hopf subalgebras.
However, we do not know of any more meaningful examples. We will
still make some statements below for Hopf algebras satisfying the
power rule, instead of just for cocommutative or commutative ones,
simply because there is no additional difficulty. We also note
that some proofs require only the power rule for positive
exponents, but we have no example where this is satisfied, but the
power rule for negative exponents is not.

\begin{Rem}\nmlabel{Remark}{powbem}
  \begin{enumerate}
    \item Since $h\pow{-1}=S(h)$, the power rule \eqref{powrule}
    specializes to $h\pow{-n}=S(h)\pow n$ for all $n\in\Z$
    when we set $m=-1$. If we also specialize $n=-1$, we get
    $S^2(h)=h$, so a Hopf algebra satisfying the power rule is
    involutive. This should be read with a grain of salt, however:
    We have already remarked that the definition of Hopf powers that we use
    only fits with the definition of exponent given by
    Etingof and Gelaki if $H$ is involutive to begin with.
    \item Conversely, it is easy to check that the power
    rule for positive exponents together with the requirement that
    $h\pow{-n}=S(h\pow n)$ for all $n\in\N$ and $S^2=\id$ implies the
    power rule for general exponents.
    \item Let $H$ be a Hopf algebra
    satisfying the power rule. If the $n$-th Hopf  power of $h\in H$
    is trivial, and $m$ is a multiple of $n$, then the $m$-th Hopf
    power of $h$ is trivial.
\end{enumerate}
\end{Rem}
\begin{Rem}
  If $H$ is finite-dimensional and satisfies the power rule, hence is
  involutive, then
  it has finite exponent. This follows from more general and
  rather deeper results of Etingof and Gelaki. In the special case
  where $H$ is commutative (or cocommutative) it was shown by
  Gabriel \cite{Gab:EISGGF}. %
  If $\ch k=0$, then by results of Larson and Radford
  \cite{LarRad:FDCHAC0S} $H$ is semisimple and cosemisimple, and
  Etingof and Gelaki \cite[Thm.4.3]{EtiGel:EFDHA} show that $H$ has
  finite exponent, which in fact divides $\dim(H)^3$. On the other
  hand, if $k$ has finite characteristic, then another result of
  Etingof and Gelaki \cite[Cor.4.10]{EtiGel:EFDHA} shows that $H$
  has finite exponent (in the same sense as in our paper, since
  $H$ is involutive).
\end{Rem}
\begin{Prop}\nmlabel{Proposition}{powruleProp}
  Let $H$ be a Hopf algebra with finite exponent $e=\exp(H)$ which
  satisfies the power rule. Then
\begin{enumerate}
  \item $\tpd_n(H)=1$ if and only if $\gcd(n,e)=1$.
  \item There is $h\in H$ of Hopf order $n$ if and only if $n|e$.
  \item \label{powruleorders}$\TPS_{e-n}(H)=\TPS_n(H)$ for $0<n<e$.
  \item \label{powrulesym}If $H$ is finite-dimensional, then
  $\tpd_{e-n}(H)=\tpd_n(H)$, and
  $\tpd_{m,e-n}(H)=\tpd_{m,n}(H)=\tpd_{e-m,e-n}(H)$ for
  $0<m,n<e$.
\end{enumerate}
\end{Prop}
\begin{proof}
Let $h\in H$. We show that if the $n$-th Hopf power of $h$ is
trivial, then so is the $k$-th Hopf power, where $k=(e,n)$ is the
greatest common divisor of $n$ and $e=\exp(H)$. Write $k=na+eb$
for some $a,b\in\Z$. Then $h\pow
k=h\pow{na+eb}=h\pow{na}=\epsilon(h)1.$ This shows that
$\TPS_n(H)$ can only contain nontrivial elements if $(n,e)\neq 1$,
and that the Hopf order of $h$ divides $e$.

Now let $n|e$. We first show that there is a nontrivial element
whose $n$-th Hopf power is trivial. Assume otherwise, that is,
assume that for every $h\in H$, $h\pow n=\epsilon(h)1$ implies
$h=\epsilon(h)1$. Then $\epsilon(h)1=h\pow e=(h\pow{e/n})\pow n$
implies that $h\pow{e/n}=\epsilon(h)1$ for all $h\in H$,
contradicting the definition of the exponent. Now assume that
there is no element of Hopf order $n$. This means that there is
$m<n$ with $\TPS_n(H)\subset\TPS_m(H)$. But then
$h\pow{e/n}\in\TPS_m(H)$ for all $h\in H$, so
$\epsilon(h)1=(h\pow{e/n})\pow m=h\pow{em/n}$ for all $h\in H$,
contradicting once more the definition of the exponent. We have
shown the other implication in (2). The missing implication in (1)
follows.

For (3) let $h\in\TPS_{e-n}$. Then $\epsilon(h)1=h\pow
n=h\pow{n-e}=S(h\pow{e-n})$ and hence $h\in\TPS_{e-n}$. Clearly
(4) follows from (3).
\end{proof}

For a more detailed comparison to the group case we make the
following definition:
\begin{Def}
  Let $G$ be a finite group. The \emph{$n$-th trivial power set} of $G$ is $$\TPM_n:=\TPM_n(G)=
  \{g\in G|g^n=\neut\}.$$ The \emph{$n$-th trivial power number} of $G$ is $\tpn_n=|\TPS_n|$.
\end{Def}
\begin{Rem}\nmlabel{Remark}{tpnrem}
  Let $G$ be a finite group with exponent $e=\exp(G)$, and $n\in\N$.
  \begin{enumerate}
    \item $\TPM_n=\{\neut\}$ if and only if $(n,e)=1$.
    \item $G$ contains an element of order $n$ if and only if
  $\TPM_n\not\subset\bigcup_{1<m<n}\TPM_m$.
    \item $\TPM_n=\TPM_{e-n}$, where $e=\exp(G)$. In particular $\tpn_n=\tpn_{e-n}$.
    \item $\TPM_n\cap\TPM_m=\TPM_{\gcd(m,n)}$ for $m,n\in\N$. In particular, $g^m=\neut$
    if and only if the order of $g$ divides $m$.\label{intersec}
  \end{enumerate}
\end{Rem}
The first property is exactly parallel to the first property in
\nmref{powruleProp}, though we will elaborate some on the
difference between $\TPS_n(kG)$ and $\TPM_n(G)$ below. The second
statement in \nmref{powruleProp} is of course simply false for
orders of group elements. The third property, which, again,
parallels \nmref{powruleProp}, perhaps deserves a short proof for
comparison: If $g^n=\neut$, then $g^{e-n}=g^e(g^n)\inv=\neut$. The
fourth property is (standard and) easy to check just using the
standard rules for taking powers in a group (write
$\gcd(m,n)=rm+sn$, and conclude that for $g\in\TPM_m\cap\TPM_n$ we
have $g^{\gcd(m,n)}=g^{rm+sn}=(g^m)^r(g^n)^s=\neut$), but its
obvious Hopf analog does not hold even for group algebras, as we
will see below.
\begin{Expl}\nmlabel{Example}{kGex}
  Let $G$ be a finite group. Then $\TPM_n(G)\subset\TPS_n(kG)$.
  Also, if $x,y\in G$ satisfy $x^n=y^n$, then $x-y\in\TPS_n(kG)$.
  More generally, consider $h=\sum_{x\in G}\alpha_xx\in kG$.
  Then $h\pow n=\sum_{x\in G}\alpha_xx^n=\sum_{g\in
  G}\left(\sum_{{\substack{x\in
  G\\x^n=g}}}\alpha_x\right)g$, and thus $h\in\TPS_n$ if and only
  if the coefficients $\alpha_x\in k$
  satisfy the equations $\sum_{{\substack{x\in
  G\\x^n=g}}}\alpha_x=0$ for all $g\in G\setminus\{\neut\}$ (or
  for all $g\in G\setminus \{\neut\}$ that are $n$-th powers to
  begin with). In particular, the space $\TPS_n(kG)$ is spanned by
  $\TPM_n(G)$ along with the differences $x-y$ with $x^n=y^n$, or
  by these differences along with the neutral element $\neut$.
  For any group $G$, the integral $\sum_{g\in G}g$ in the Hopf
  algebra $kG$ has Hopf order $\exp(kG)$.
\end{Expl}
\begin{Expl}\nmlabel{Example}{gralgexex}
\begin{enumerate}
  \item In the symmetric group $S_4$, the elements $(1234)$ and
  $(1432)$, each of order $4$, have the same square $(13)(24)$, so
  their difference has Hopf order $2$. The same holds for
  $(1243)-(1342)$ and $(1423)-(1324)$. It is not hard to check that
  $\dim\TPS_2(kS_4)=|\TPM_2(S_4)|+3$.\label{moretriv}
  \item In the symmetric group $S_5$, consider the elements
  \begin{align*}
    x&=(12)(345)&y&=(345)\\
    a&=(12)(354)&b&=(354).
  \end{align*}
  Since $x^2=y^2=(354)$, and $a^2=b^2=(345)$, the Hopf order of
  $h:=x-y+b-a\in kS_5$ is $2$. On the other hand $x^3=a^3=(12)$ and
  $y^3=b^3=\id$, and thus also $h\pow
  3=0=\epsilon(h)1$. In particular $\tpd_{2,3}(kS_5)$\label{twoandthree}
\end{enumerate}
\end{Expl}
\begin{Expl}\nmlabel{Example}{kG*ex}
  Let $G$ be a finite group. Consider the Hopf algebra $H=k^G$ of $k$-valued functions on $G$, the dual of the group algebra $kG$.
  Let $\bas g$ for $g\in G$ denote the elements of the canonical
  basis of $k^G$.
\begin{enumerate}
  \item\label{kG*triv} Since $h\sw 1(x)h\sw 2(y)=h(xy)$ and $(hh')(x)=h(x)h'(x)$ for $h,h'\in H$ and $x,y\in G$,
  we have $h\pow n(x)=(h\sw 1\cdot\dots\cdot h\sw n)(x)=h\sw 1(x)\cdot\dots\cdot h\sw n(x)=h(x^n)$.
  Since on the other hand $(\epsilon(h)1)(x)=h(\neut)$, we see that the $n$-th Hopf power of $h\in k^G$
  is trivial if and only if $h(x^n)=h(e)$ for all $x\in G$. In other words, $\TPS_n(k^G)$
  consists precisely of those functions that are
  constant on the set of all $n$-th powers in $G$.
  \item Let $\neut\neq g\in G$ and $n\in\N$. Then $(\bas g)\pow n=0=\epsilon(\bas g)1$ if and only
  if $g$ is not an $n$-th power in $G$. In particular, the Hopf order of $\bas g$ is the least
  positive integer such that $g$ is not an $n$-th power in $G$.
  \item\label{kG*int} The $n$-th Hopf power of $\bas\neut$ is trivial if and only if no element besides $\neut$
  is an $n$-th power in $G$, in other words, if and only if all $n$-th powers equal $\neut$.
  Thus the Hopf order of the integral  $\bas \neut$ of $k^G$ is the exponent of $G$.
  \item Let $p$ be a prime divisor of the order of $G$. Then there is $g\in G$ such that
  $\bas g$ has Hopf order $p$. More generally, if $G$ contains an element of order $p^k$,
  then there is $g\in G$ such that $\bas g$ has Hopf order $p^k$.
\end{enumerate}
\end{Expl}
\begin{proof}
We need only prove the last assertion.

Note that in a finite abelian group $A$, every element of $A$ is
an $m$-th power if $m$ is prime to the order of $A$. Indeed, in
the abelian case the map $A\ni a\mapsto a^m\in A$ is a group
homomorphism, and if $g$ is in its Kernel, that is, $g^m=e$, then
$m$ is a multiple of the order of $g$, which divides the order of
$A$.

Now choose $g\in G$ of order $p^k$, where $k$ is maximal. Then
$\bas g$ has Hopf order $p$. In fact $g$ is not a $p$-th power by
maximality of $k$. On the other hand, any $m<p$ is prime to $p$,
and so $g$ is an $m$-th power even in the subgroup of order $p^k$
generated by $g$.

More generally, for any $1\leq n\leq k$, the element
$g^{p^{n-1}}\in G$ has order $p^{k-n+1}$, hence is not a $p^n$-th
power in $G$ by maximality of $k$. Let $1<m<p^n$. If $m$ is a
power of $p$ then $g^{p^{n-1}}$ is obviously an $m$-th power of a
power of $g$. Otherwise write $m=p^t\ell$ with $\ell$ prime to $p$
and $0\leq t<n$. Then $g$ is an $\ell$-th power in the subgroup of
order $p^k$ generated by $g$, say $g=h^\ell$,  and hence
$g^{p^{n-1}}=h^{\ell p^{n-1}}=(h^{p^{n-1-t}})^m$.
\end{proof}

Before discussing a few examples, we will make some observations
on general finite-dimensional Hopf algebras. From the descriptions
of $\TPS_n(k^G)$ and $\TPS_n(kG)$ obtained in \nmref{kGex} and
\nmref{kG*ex} \eqref{kG*triv}, it is easy to check that these
spaces have the same dimension. More generally, $\TPS_n(H)$ and
$\TPS_n(H^*)$ have the same dimension for any finite-dimensional
Hopf algebra $H$, since these spaces are the kernels of an
endomorphism of $H$ and its dual map. This is \eqref{tdual} of the
following \ref{forduals.nme}, and also follows directly from
\eqref{ttp}, which is based on a closer analysis of the spaces
$\TPS_n(H)$.
\begin{Lem}\nmlabel{Lemma}{forduals}
  Let $H$ be a finite-dimensional Hopf algebra.
  \begin{enumerate}
    \item $\TPS_n(H)=k\cdot 1_H\oplus\Ker([n])$.\label{tstruc}
    \item $\tpd_n(H)=\dim(H)+1-\rank([n])$.
    \item $\tpd_n(H)=\tpd_n(H^*)$.\label{tdual}
    \item If $K$ is another finite-dimensional Hopf algebra,
    then\label{ttp}
    \begin{imini}
      \begin{multline*}\tpd_n(H\o K)=\tpd_n(H\o
      K^*)\\=(\tpd_n(H)-1)\dim(K)+(\tpd_n(K)-1)\dim(H)+1-(\tpd_n(H)-1)(\tpd_n(K)-1)
      \end{multline*}
    \end{imini}
  \end{enumerate}
\end{Lem}
\begin{proof}
    Since $[n](1)=1$ and $\epsilon[n]=\epsilon$, the $n$-th power
  map preserves the direct sum decomposition $H=k\cdot 1_H\oplus
  \Ker(\epsilon)$, and it is the identity on the first summand.
  This implies (1) and (2). Since $[n]_{H\o K}=[n]_H\o[n]_K$, and
  $[n]_{K^*}=([n]_K)^*$, the formula in (4) follows by substituting the rank formula
  $\rank(f\o g)=\rank(f)\rank(g)$ for homomorphisms $f,g$.
\end{proof}
\begin{Rem}
  In the examples we have computed, not only
  $\tpd_n(kG)=\tpd_n(k^G)$ and $\tpd_n(kG\o kF)=\tpd_n(k^G\o kF)$
  for finite groups $F,G$, which is proved in greater generality
  in the preceding Lemma, but also
  $\tpd_{m,n}(kG)=\tpd_{m,n}(k^G)$ and $\tpd_{m,n}(kG\o
  kF)=\tpd_{m,n}(k^G\o kF)$. We do not know if this holds for all
  finite groups $F,G$.
\end{Rem}
We will now list the numbers $\tpd_{i,j}$ for a few group
algebras, and the trivial power numbers of the respective group
for comparison. We begin with the symmetric group $S_3$ and the
alternating group $A_4$, and discuss at the same time the format
of tables we are also planning to use for more general examples of
Hopf algebras. At the bottom of each of the tables
$(\tpd_{i,j}(\QQ G))$ we have attached a list of the numbers
$\tpn_i(G)$. This of course only makes sense for group algebras.
\begin{table}%[h]
\label{twogps}

\begin{align*}
 \begin {array}{|c|ccccc|}
 \hline \tpd_{i,j}(\QQ S_3)&j=1&2&3&4&5\\\hline1&1&1&1&1&1\\2&
 &4&1&4&1\\3& & &3 &1&1\\4& & & &4&1\\5& & & & &1\\\hline
 \noalign{\bigskip}\hline \tpn_i(S_3)&1&4&3&4&1\\\hline\end
 {array}
&&&
\begin {array}{|c|ccccc|}
 \hline
 \tpd_{i,j}(\QQ A_4)&j=1&2&3&4&5\\
 \hline
 i=1&1&1&1&1&1\\
 2&&4&1&4&1\\
 3& & &9 &1&1\\
 4& & & &4&1\\
 5& & & & &1\\
 \hline
 \noalign{\bigskip}
 \hline
 \tpn_i(A_4)&1&4&9&4&1\\
 \hline
\end {array}
\end{align*}
\caption{}
\end{table} We have only printed the upper triangular part of the tables
$(\tpd_{i,j}(H))$; for any Hopf algebra $H$ the obvious symmetry
$\tpd_{m,n}(H)=\tpd_{n,m}(H)$ makes the lower part redundant. Any
rows or columns beyond the $e$-th, where $e$ is the exponent (in
this case, $e=6$ for both groups) would be redundant since
$\TPS_{n+e}(H)=\TPS_n(H)$ by \nmref{genpowrem}. The $e$-th row and
column are also redundant because $\TPS_e(H)=H$ and hence
$\tpd_{e,n}(H)=\tpd_n(H)=\tpd_{n,n}(H)$; that is, the $e$-th
column as well as the diagonal contains the numbers $\tpd_n(H)$.
We have kept the first row and $(e-1)$-st column, which will
contain only $1$'s for any finite-dimensional Hopf algebra: The
first row because the first Hopf power is the identity, the last
column because the $(e-1)$-st Hopf power map is the antipode. For
group algebras we have the additional symmetries in
\nmref{powruleProp}, \eqref{powrulesym} that would allow us to
further reduce the table; one reason we didn't do this is that
only some of these symmetries will continue to hold in interesting
more general examples. The other is that the table in the above
form gives a convenient recipe to look up the numbers $n$ for
which elements of Hopf order $n$ exist. We will return to this
below.

Note that in the above examples $\tpn_i(G)=\tpd_i(\QQ G)$. We
already know from \nmref{gralgexex} \eqref{moretriv} that this
will fail for $G=S_4$. Table~\ref{S4group} lists the numbers for
this case.
\begin{table}%[h]
\label{S4group}
$$\begin{array}{|c|ccccccccccc|} \hline
\tpd_{i,j}(\QQ S_4)&j=1&2&3&4&5&6&7&8&9&10&11\\ \hline
i=1&1&1&1&1&1&1&1&1&1&1&1
\\ 2&&13&1&13&1&13&1&13&1&13&1\\
3& &&9&1&1&9&1&1&9&1&1
\\
4&&&&16&1&13&1&16&1&13&1
\\
5&&&&&1&1&1&1&1&1&1\\
6&&&&&&21&1&13&9&13&1\\
7&&&&&&&1&1&1&1&1\\
8&&&&&&&&16&1&13&1\\
9&&&&&&&&&9&1&1\\
10&&&&&&&&&&13&1\\
11&&&&&&&&&&&1\\ \hline \noalign{\bigskip}
 \hline
\tpn_i(S_4)&1&10&9&16&1&18&1&16&9&10&1\\\hline
\end {array}$$
\caption{}
\end{table}
 For each of the groups $G=S_3,A_4,S_4$ we see that
$\tpd_{i,j}(kG)=\tpd_{\gcd(i,j)}(kG)$ for all $i,j$. This looks
like a variant of the equation
$|\TPM_i(G)\cap\TPM_j(G)|=\tpn_{\gcd(i,j)}(G)$ that follows from
\nmref{tpnrem} \eqref{intersec}, and is the reason that we only
print the numbers $\tpn_i(G)$ for a group, rather than a full
square table of numbers. The observation on group algebras is
misleading, however, as we have seen in \nmref{gralgexex}
\eqref{twoandthree}.

As we mentioned above, it is easy to read off from the tables for
which numbers $n$ elements of Hopf order $n$ exist. We will spell
out the recipe, which rephrases \nmref{ordersbydim}, for later
reference, and we will illustrate it with the example $kS_4$,
although we know the result here without looking from
\nmref{powruleProp} \eqref{powruleorders}.
\begin{Rem}\nmlabel{Remark}{recipe}
  For a finite dimensional Hopf algebra $H$ with finite exponent
  $e$ consider the upper triangular table $(\tpd_{i,j}(H))$ for
  $1\leq j\leq i<e$. The Hopf algebra $H$ contains a nontrivial
  element whose $n$-th power is trivial if and only if the $n$-th
  diagonal entry is greater than $1$. It contains an element of
  Hopf order $n$ if and only if that diagonal entry is strictly
  larger than all the entries above it.
\end{Rem}
For example, we see that $\QQ S_4$ contains nontrivial elements
whose tenth Hopf power is trivial, since in Table~\ref{S4group}
the tenth diagonal entry is $13$. But there is no element of Hopf
order $10$, since the same number $13$ appears several times above
the diagonal in the tenth column. By contrast, there do exist
elements of Hopf order six, since the sixth diagonal element $21$
is strictly larger than all the numbers in the column above it.

\section{Doubles and twists}\nmlabel{Section}{sec:doubles}

Let $H$ be a bialgebra, and $\theta\in H\o H$ an invertible
element. If $\theta$ satisfies the cocycle identity
\begin{equation*}
  (\theta\o 1)\cdot(\Delta\o
  H)(\theta)=(1\o\theta)\cdot(H\o\Delta(\theta),
\end{equation*}
then $\Delta_\theta(h)=\theta\Delta(h)\theta\inv$ defines a new
comultiplication making the algebra $H$ into a bialgebra
$H_\theta$, called a \emph{Drinfeld twist} \cite{Dri:QHA} of $H$.
If $(\epsilon\o H)(\theta)=(H\o\epsilon)(\theta)=1$, the counit of
$H_\theta$ is that of $H$. If $H$ is a Hopf algebra, so is
$H_\theta$.

The dual construction was studied by Doi \cite{Doi:BBQB}: If a
convolution invertible map $\sigma\colon H\o H\to k$ is a
\emph{two-cocycle}, that is, satisfies the identity
\begin{equation*}
  \sigma(f\sw 1\o g\sw 1)\sigma(f\sw 2g\sw 2\o h)=\sigma(g\sw 1\o
  h\sw 1)\sigma(f\o g\sw 2h\sw 2)
\end{equation*}
for all $f,g,h\in H$, then
$$g\cdot h:=\sigma(g\sw 1\o h\sw 1)g\sw 2h\sw 2\sigma\inv(g\sw 3\o h\sw 3)$$
defines a new multiplication making the coalgebra $H$ into a
bialgebra $H^\sigma$, which we call a cocycle twist of $H$. If
$\sigma(1\o h)=\sigma(h\o 1)=\epsilon(h)$ for all $h$, then $1_H$
is the unit of $H^\sigma$. If $H$ is a Hopf algebra, so is
$H^\sigma$.

Etingof and Gelaki \cite{EtiGel:EFDHA} show that their improved
version of exponent of a finite-dimensional Hopf algebra is
invariant under Drinfeld twists (and cocycle twists, since it is
invariant under taking the dual). They also show that the exponent
is invariant under the operation of taking the Drinfeld double;
this is actually a special case of the result on cocycle twists:
Doi and Takeuchi \cite{DoiTak:MATCQV} have shown that the Drinfeld
double can be constructed in a particularly smooth way as a
cocycle twist. In this section we will start to investigate
whether the trivial power dimensions $\tpd_n(H)$, the property of
a Hopf algebra to contain an element of Hopf order $n$, or the
property to contain nontrivial elements whose $n$-th power is
trivial, are invariant under Drinfeld twists, more particularly
under those twists that construct Drinfeld doubles.

Let us recall first from \cite{DoiTak:MATCQV} how certain cocycle
twists including the double arise from skew pairings: A skew
pairing between bialgebras $B,H$ is a map $\tau\colon B\o H\to k$
that satisfies
\begin{align*}
  \tau(bc\o h)&=\tau(b\o h\sw 1)\tau(c\o h\sw 2)\\
  \tau(b\o gh)&=\tau(b\sw 1\o h)\tau(b\sw 2\o g)\\
  \tau(b\o 1)&=\epsilon(b)\\
  \tau(1\o h)&=\epsilon(h).
\end{align*}
If $\tau$ is convolution invertible, its inverse satisfies
\begin{align*}
  \tau\inv(bc\o h)&=\tau(c\o h\sw 1)\tau(b\o h\sw 2)\\
  \tau\inv(b\o gh)&=\tau(b\sw 1\o g)\tau(b\sw 2\o h)\\
  \tau\inv(b\o 1)&=\epsilon(b)\\
  \tau\inv(1\o h)&=\epsilon(h).
\end{align*}
For an invertible skew pairing $\tau$, the map
$$[\tau]\colon B\o H\o B\o H\ni b\o g\o c\o h\mapsto
\epsilon(b)\epsilon(h)\tau(c\o g)\in k$$
is a two-cocycle on the tensor product bialgebra $B\o H$.
We denote the cocycle twist by
$$B\bowtie_\tau H:=(B\o H)^{[\tau]}$$
and note
$$(b\bowtie g)(c\bowtie h)=b\tau(c\sw 1\o g\sw 1)c\sw 2\bowtie
g\sw 2\tau\inv(c\sw 3\o g\sw 3)h.$$

If $H$ is a finite-dimensional Hopf algebra, there is an obvious
skew pairing $\tau\colon (H^*)^\cop\o H\to k$. The Drinfeld double
of $H$ can be obtained as
$$D(H)=(H^*)^\cop\bowtie_{\tau}H.$$

Table~\ref{DS3}
\begin{table}[h]
$$
\begin {array}{|c|ccccc|}
  \hline
  \tpd_{i,j}(D(\QQ S_3))&j=1&2&3&4&5\\
  \hline
  i=1&1&1&1&1&1\\2&&25&13&23&
1\\3&&&21&13&1\\4&&&&25&1
\\5&&&&&1\\\hline\end {array}
$$\caption{}\label{DS3}
\end{table}
lists the results of our computations of the
numbers $\tpd_{m,n}(H)$, where $H$ is the Drinfeld double $D(\QQ
S_3)$. Comparing with Table~\ref{S3S3}, which
\begin{table}[hb]
$$\begin {array}{|c|ccccc|}
 \hline
\tpd_{i,j}(\QQ^{S_3}\o \QQ S_3)&j=1&2&3&4&5\\
\hline
  i=1&1&1&1&1&1\\2&&28&13&28&1\\3&&&21&13&1\\4&&&&28&1
\\5&&&&&1\\\hline\end {array}
$$
\caption{}\label{S3S3}
\end{table}
does the same for
the tensor product $H=\QQ^{S_3}\o\QQ S_3$, we see that the numbers
$\tpd_n$ are not invariant under the specific cocycle twist that
obtains the double from the tensor product.

Our calculations have shown that $\tpd_{i,j}(H)=\tpd_{i,j}(H^*)$
for $H=D(\QQ S_3)$, so the analogous table for $D(\QQ S_3)^*$
looks exactly like Table~\ref{DS3}. In particular, the numbers
$\tpd_{i}(H)$ are not invariant under Drinfeld twists either.

Applying \nmref{recipe} we see that $D(\QQ S_3)$ has elements of
Hopf order four. By \nmref{powruleProp} this is impossible for a
Hopf algebra of exponent six satisfying the power rule --- such as
$\QQ^{S_3}\o \QQ S_3$, from which $D(\QQ S_3)$ is obtained by a
cocycle twist.

Let us also point out the inequality $\tpd_{2,4}(D(\QQ
S_3))=23<25=\tpd_2(D(\QQ S_3))$, which shows that $D(\QQ S_3)$
contains elements whose second Hopf power is trivial, but whose
fourth Hopf power is not.

Our next set of examples is based on the alternating group $A_4$.
\begin{table}[h]
$$
\begin {array}{|c|ccccc|}
 \hline \tpd_{i,j}(\QQ^{A_4}\o \QQ A_4)&j=1&2&3&4&5\\
 \hline
 i=1&1&1&1&1&1\\
 2&&64&49&64&1\\
 3&&&129&49&1\\
 4&&&&64&1\\
 5&&&&&1\\
 \hline\end {array}
$$\caption{}\label{A4A4}
\end{table}
Table~\ref{A4A4} lists the numbers $\tpd_{i,j}(\QQ ^{A_4}\o\QQ
A_4)$. As in the example of $S_3$, the numbers for the cocycle
twist $D(\QQ A_4)$, listed in
\begin{table}[h]
$$
\begin {array}{|c|ccccc|}
 \hline
 \tpd_{i,j}(D(\QQ A_4))&j=1&2&3&4&5\\
 \hline
 i=1&1&1&1&1&1\\2&&64&49&34
 &1\\3&&&121&49&1\\4&&&&64&1
 \\5&&&&&1\\\hline\end {array}
$$\caption{}\label{DA4}
\end{table}
Table~\ref{DA4}, are different. From
\begin{table}[h]
$$
\begin {array}{|c|ccccc|}
 \hline
 \tpd_{i,j}(D(\QQ A_4)^*)&j=1&2&3&4&5\\
 \hline
 i=1&1&1&1&1&1\\2&&64&49&40
 &1\\3&&&121&49&1\\4&&&&64&1
 \\5&&&&&1\\\hline\end {array}
$$\caption{}\label{dDA4}
\end{table}
Table~\ref{dDA4} we see that $\tpd_{i,j}(H)\neq\tpd_{i,j}(H^*)$
for $H=D(\QQ A_4)$. Observe, however, that the Hopf orders of
elements that occur are still the same for $H$ and $H^*$ here: In
each case all Hopf orders in $\{1,...,4\}$ are possible.

Each of the tables we have seen so far is symmetric with respect
to the anti-diagonal; this means we have
$\tpd_{e-m,e-n}=\tpd_{n,m}$, or equivalently
$\tpd_{e-m,e-n}=\tpd_{m,n}$ for all $m,n$. For Hopf algebras that
satisfy the power rule this follows from the more general rule
$\tpd_{m,n}=\tpd_{m,e-n}$; both were observed in
\nmref{powruleProp} \eqref{powrulesym}. The more general symmetry
$\tpd_{m,n}=\tpd_{m,e-n}$  does not hold for doubles of group
algebras: For example $\tpd_{2,2}(D(\QQ A_4))\neq\tpd_{2,4}(D(\QQ
A_4))$ although $4=e-2$. We will show, however, that the observed
mirror symmetry about the anti-diagonal does hold for a large
class of Hopf algebras. Of course once we prove that, part of the
tables is actually redundant: We could have left out the part
below, or the part above the anti-diagonal. We will make use of
this for bigger tables below, but whenever possible we have left
the redundant data in the picture to facilitate the procedure of
looking up the numbers $n$ for which elements of Hopf order $n$
exist.

\begin{Def}
  Let $k$ be a field with involution $k\ni x\mapsto \ol x\in k$. A
  \emph{$*$-Hopf algebra} over $k$ is a Hopf algebra $H$ with a
  semilinear involution $*\colon H\to H$ on the algebra $H$ such that
  $\Delta\colon H\to H\o H$ and $\epsilon\colon H\to k$ are
  homomorphisms of $*$-algebras. Semisimple $*$-Hopf algebras over
  the complex numbers are called \emph{Kac Hopf algebras}.
\end{Def}
\begin{Prop}\nmlabel{Proposition}{symmetry}
  Let $H$ be a finite-dimensional
  Hopf algebra with finite exponent $e$, and
  $1\leq n<e$.
  \begin{enumerate}
    \item $\TPS_n(H)=\TPS_{e-n}(H^\op)$.\label{TSHHop}
    \item $\exp(H^\op)=\exp(H)$.\label{eHHop}
    \item If $H$ is a $*$-Hopf algebra, then $\TPS_n(H)$ and
    $\TPS_{e-n}(H)$ correspond to each other under the involution
    of $H$.\label{forKac}
    \item If $H$ is a $*$-Hopf algebra, or isomorphic to its
    opposite Hopf algebra, then $\tpd_n(H)=\tpd_{e-n}(H)$ and
    $\tpd_{m,n}(H)=\tpd_{e-m,e-n}(H)$ for $1\leq
    m,n<e$.\label{dimsym}
  \end{enumerate}
\end{Prop}
\begin{proof}
  Let $0\leq n\leq e$, and $h\in \TPS_n(H)$. Then $\epsilon(h)1=h\pow
  n=h\pow{n-e}=S(h\sw 1)\cdot\dots\cdot S(h\sw {e-n})$. Applying
  $S\inv$ yields $\epsilon(h)1=h\sw {e-n}\cdot\dots\cdot h\sw 1$,
  and thus $\TPS_n(H)\subset \TPS_{e-n}(H^\op)$. In the
  special case $n=0$ this means $\exp(H^\op)\leq e= \exp(H)$.
  By symmetry this proves \eqref{TSHHop} and \eqref{eHHop}.

  Now if we have a Hopf algebra isomorphism $f\colon H\to H^\op$, then
  \eqref{TSHHop} implies $f(\TPS_n(H))=\TPS_{e-n}(H)$, and hence
  the claims in \eqref{dimsym}. If \eqref{forKac} and the claims
  in \eqref{dimsym} for $*$-Hopf algebras are not strictly a
  special case of this, it is only for a technical reason: The
  involution of $H$ is not linear. Still, if we apply the
  involution to $\epsilon(h)1=h\sw{e-n}\cdot\dots\cdot h\sw 1$ for
  $h\in\TPS_n(H)=\TPS_{e-n}(H^\op)$, we get
  $\epsilon(h^*)1=\ol{\epsilon(h)}1=\left(h\sw{e-n}\cdot\dots\cdot
  h\sw 1\right)^*=(h\sw 1)^*\cdot\dots\cdot (h\sw{e-n})^*=(h^*)\sw
  1\cdot\dots\cdot(h^*)\sw {e-n}$, so again
  $(\TPS_n(H))^*\subset\TPS_{e-n}(H)$, and equality follows by
  symmetry.
\end{proof}
The equality in \eqref{eHHop} was obtained for the improved
version of exponent by Etingof and Gelaki
\cite[Cor.2.6]{EtiGel:EFDHA}, but with a less elementary proof. It
is also proved in a different way by Kashina, Sommerh\"auser, and
Zhu \cite{KasSomZhu:HFSI}.

Since the Drinfeld doubles of finite groups over $\mathbb C$ are
Kac algebras, \nmref{symmetry} \eqref{forKac} explains the mirror
symmetry we observed in the tables above.
\begin{table}
$$\begin {array}{|c|ccccccccccc|}
 \hline \tpd_{i,j}(D(\QQ S_4))&1&2&3&4&5&6&7&8&9&10&11
\\\hline1&1&1&1&1&1&1&1&1&1&1&1\\2& &
415&217&395&1&401&1&401&217&386&1\\3& & &313&
241&1&305&1&241&304&217&1\\4& & & &484&1& 449&1&443&241&401&1\\5&
& & & &1&1&1&1&1 &1&1\\6& & & & & &535&1&449&305&401&1
\\7& & & & & & &1&1&1&1&1
\\8& & & & & & & &484&241&395&1
\\9& & & & & & & & &313&217&
1\\10& & & & & & & & & & 415&1\\11& & & & & & & & &  &
&1\\\hline\end {array}$$

\caption{}\label{S4dbl}
\end{table}

Our next example is the Drinfeld double of the symmetric group
$S_4$. As in the previous examples, we can spot many Hopf orders
from Table~\ref{S4dbl} which would be impossible for Hopf algebras
that satisfy the power rule. However, the primes $5$ and $7$ are
not possible Hopf orders; there do not even exist nontrivial
elements whose fifth or seventh Hopf powers are trivial. Thus, for
the group $G=S_4$, we see that nontrivial elements of $D(\QQ G)$
whose $n$-th Hopf power is trivial exist only for those $n$ where
the tensor product $\QQ^{G}\o\QQ G$, or, for that matter, the
group algebra $\QQ G$, already contains such elements.

This limited version of an ``invariance'' of the behavior of
trivial power dimensions under special cocycle twists is proved
for arbitrary groups in the following \ref{grdbl.nme}. It rules
out, in particular, the existence of elements of a double $D(kG)$
whose Hopf order is a prime that does not divide the dimension.

\begin{Prop}\nmlabel{Proposition}{grdbl}
  Let $\tau\colon B\o H\to k$ be an invertible skew pairing of
  bialgebras.

  If $n\in\N$ satisfies $\tpd_n(H)=\tpd_n(B)=1$, and $H$ is cocommutative,
  then also $\tpd_n(B\bowtie_\tau H)=1$.

  In particular, if $G$ is a finite group, and $n$ is prime to the order of $G$,
  then the only elements in the double $D(kG)$ that have trivial $n$-th Hopf power
  are the scalar multiples of the identity.
\end{Prop}
\begin{proof}
  We first essentially follow Kashina \cite[Sec.~3]{Kas:GPMHA}
  to find a formula for the $n$-th power
  map on $B\bowtie H$. We rewrite multiplication in $B\bowtie H$
  in the form
  \begin{align*}
    (b\bowtie g)(c\bowtie h)&=b\tau(c\sw 1\o g\sw 1)c\sw 2\bowtie g\sw
    2\tau\inv(c\sw 3\o g\sw 3)h\\
    &=b(g\sw 1\hit c\sw 1)\bowtie(g\sw 2\hitby c\sw 2)h
  \end{align*}
  Using the actions
  \begin{align*}
     H\o B\ni h\o b\mapsto h\hit b&=\tau(b\sw 1\o h)b\sw 2\in
     B\\
     H\o B\ni h\o b\mapsto h\hitby b&=\tau\inv(b\o h\sw 2)h\sw 1\in
     H.
  \end{align*}
  It is easy to check that $\hit$ makes $B$ an $H$-module algebra,
  satisfying in addition $\Delta(h\hit b)=h\hit b\sw 1\o b\sw 2$,
  and that $\hitby$ makes $H$ a $B$-module algebra, satisfying in
  addition $\Delta(h\hitby b)=h\sw 1\o h\sw 2\hitby b$.
  We consider the bijection
  \begin{align*}
   F\colon B\o H &\to B\bowtie H\\
   b\o h&\mapsto \tau(b\sw 1\o h\sw 1)b\sw 2\bowtie h\sw
2\tau\inv(b\sw 3\o h\sw 3)\\&=h\sw 1\hit b\sw 2\bowtie h\sw
2\hitby b\sw 2\\&=(1\bowtie h)(b\bowtie 1)
  \end{align*}
  and will prove, by induction on $n$,
  $$[n]F(b\o h)=h\sw 1\hit b\sw 1\pow n\bowtie h\sw 2\pow
  n\hitby b\sw 2.$$
  The case $n=1$ is just the definition of multiplication in $B\bowtie H$.
  Assuming the formula for $n$, we first get
  \begin{align*}
  (b\bowtie h)\pow{n+1}
    &=(b\sw 1\bowtie 1)\bigl((1\bowtie h\sw 1)(b\sw 2\bowtie 1)\bigr)\pow{n}(1\bowtie h\sw 2)\\
    &=b\sw 1(h\sw 1\hit b\sw 2\pow n)\bowtie(h\sw 2\pow n\hitby b\sw 3)h\sw 3
  \end{align*}
  In this expression we substitute $h\sw 1\hit b\sw 1\bowtie h\sw 2\hitby b\sw 2$
  for $b\bowtie h$,
  using
  $$\Delta(h\sw 1\hit b\sw 1)\o\Delta(h\sw 2\hitby b\sw 2)
  =h\sw 1\hit b\sw 1\o b\sw 2\o h\sw 2\o h\sw 3\hitby b\sw 3,$$
  and find
  \begin{align*}
    [n+1]F(b\o h)&=(h\sw 1\hit b\sw 1\bowtie h\sw
    2 \hitby b\sw 2)\pow{n+1}\\
    &=(h\sw 1\hit b\sw 1)(h\sw 2\hit b\sw 2\pow n)\bowtie(
    h\sw 3\pow n\hitby b\sw 3)(h\sw 4\hitby b\sw 4)\\
    &=h\sw 1\hit b\sw 1b\sw 2\pow n\bowtie h\sw 2\pow nh\sw
    3\hitby  b\sw 3\\
    &=h\sw 1\hit b\sw 1\pow{n+1}\bowtie h\sw 2\pow{n+1}\hitby  b\sw 2
  \end{align*}
  Now consider the bijections
  \begin{align*}
    P\colon B\o H\ni b\o h&\mapsto h\sw 1\hit b\o h\sw 2\in B\o H\\
    Q\colon B\o H\ni b\o h&\mapsto b\sw 1\o h\sw 1\tau\inv(b\sw 2\o
    h\sw 2)\in B\o H.
  \end{align*}
  It is easy to check that
  $$Q^n(b\o h)= b\sw 1\o h\sw 1\tau\inv(b\sw 2\o h\sw 2\pow n).$$
  If $H$ is cocommutative, then $[n]_H$ is a coalgebra map, and we
  finally have
  \begin{align*}
    [n]_{B\bowtie H}F(b\o h)&=h\sw 1\hit b\sw 1\pow n\tau\inv(b\sw 2\o h\sw
    3\pow n)\bowtie h\sw 2\pow n\\
    &=(B\o[n]_{H})P([n]_{B}\o H)Q^n(b\o h).
  \end{align*}
  In particular, if $\tpd_n(H)=\tpd_n(B)=1$, so that $[n]_B$ and
  $[n]_H$ are bijective, then $[n]_{B\bowtie H}$ is bijective, and
  so $\tpd_n(B\bowtie H)=1$.
\end{proof}
\begin{table}
$$\begin {array}{|c|ccccccccccc|}
 \hline \tpd_{i,j}(D(\QQ S_4)^*)&1&2&3&4&5&6&7&8&9&10&11
\\\hline 1&1&1&1&1&1&1&1&1&1&1&1\\2& &
415&205&383&1&415&1&415&205&359&1\\3& & &313&
241&1&301&1&241&307&205&1\\4& & & &484&1& 463&1&452&241&415&1\\5&
& & & &1&1&1&1&1 &1&1\\6& & & & & &535&1&463&301&415&1
\\7& & & & & & &1&1&1&1&1
\\8& & & & & & & &484&241&383&1
\\9& & & & & & & & &313&205&
1\\10& & & & & & & & & & 415&1\\11& & & & & & & & & & &1
\\\hline
\end{array}$$

\caption{}\label{dS4dbl}
\end{table}

To close the section, note that according to Table~\ref{dS4dbl},
the dual $D(\QQ S_4)^*$ of the double of $S_4$ does not contain an
element of Hopf order $10$, whereas $D(\QQ S_4)$ does.
\section{Bismash products}\nmlabel{Section}{sec:bismash}

We will now compute the numbers $\tpd_{m,n}(H)$ for some bismash
products obtained from a matched pair of groups, obtained in turn
from a factorizable group. We will recall these notions and
relevant facts, using Masuoka's survey \cite{Mas:HAEC} as a
general reference.

A group $L$ is called \emph{factorizable} into subgroups
$F,G\subset L$ if $FG=L$ and $F\cap G=\{\neut\}$, that is, if
every $\ell\in L$ can be written uniquely as a product $\ell=ax$
with $a\in F$ and $x\in G$.

A \emph{matched pair of groups} $(F,G,\Hit,\Hitby)$ is a pair of
groups $F,G$ with group actions
$$G\xleftarrow{\Hitby}G\times F\xrightarrow{\Hit}F$$
of each group on the underlying set of the other, such that
\begin{align*}
  x\Hit ab&=(x\Hit a)((x\Hitby a)\Hit b\\
  xy\Hitby a&=(x\Hitby (y\Hit a))(y\Hitby a)
\end{align*}
hold for all $a,b\in F,x,y\in G$.

A factorizable group gives rise to a matched pair of groups, if we
define the mutual actions by the formula $xa=(x\Hit a)(x\Hitby a)$
for $a\in F,x\in G$, that is, let $x\Hit a\in F$ and $x\Hitby a\in
G$ be the --- by hypothesis unique --- elements of $F$ and $G$
whose product is $xa$.

Conversely, given a matched pair $(F,G,\Hit,\Hitby)$, we can
define a group $F\bowtie G$ with underlying set $F\times G$ and
multiplication $$(a,x)(b,y)=(a(x\Hit b),(x\Hitby b)y)$$ for
$a,b\in F$, $x,y\in G$. This group $L=F\bowtie G$ is factorizable
into the subgroups $F'=F\times\{\neut\}\cong F'$ and
$G'=\{\neut\}\times G\cong G'$.

Let $(F,G,\Hit,\Hitby)$ be a matched pair, with $G$ finite. Then
one can define a Hopf algebra structure on $H=k^G\o kF=:k^G\# kF$,
called a \emph{bismash product} of $k^G$ and $kF$ by
\begin{align*}
  (\bas x\#a)(\bas y\#b)&=\delta_{x\Hitby a,y}\bas x\#ab\\
  \Delta(\bas x\#a)&=\sum_{y\in G}\bas{xy\inv}\# y\Hit a\o\bas y\#
  a\\
  \epsilon(\bas x\# a)&=\delta_{\neut,x}
\end{align*}
for $a,b\in F,x,y\in G$; the unit element is $1\# \neut=\sum_{x\in
G}\bas x\#\neut$. We can view $k^G\cong k^G\# k\cdot 1$ as a Hopf
subalgebra, and $kF\cong k\cdot 1\# kF\subset k^G\# kF$ as a
subalgebra of $k^G\# kF$.

A special case of the bismash product construction is the Drinfeld
double of the group algebra of a finite group $G$. The relevant
matched pair is $(G,G,\Hit,\Hitby)$ with trivial action $\Hit$,
and adjoint action $x\Hitby a=a\inv xa$. The actions can be viewed
as coming from the factorization of the group $L=G\times G$ into
two subgroups isomorphic to $G$, the first being the diagonal
$\{(g,g)|g\in G\}$, the second $\{\neut\}\times G$. The bismash
product $k^G\# kG$ is isomorphic to the Drinfeld double $D(kG)$.

We return to the general case of a matched pair
$(F,G,\Hit,\Hitby)$. The bismash product $k^G\# kF$ is the neutral
element in the Opext group of Hopf algebra extensions
\begin{equation}
  k^G\to H\to kF.\label{extpic}
\end{equation}
 The general middle term of such a short
exact sequence is a bicrossproduct, which has its multiplication
and comultiplication deformed in addition by two $2$-cocycles. If
$k=\C$, then any middle term of an extension \eqref{extpic} is a
Kac algebra by a result of Masuoka \cite[Rem.2.4]{Mas:FFFCHAE}.
For the bismash product itself, the $*$-structure can be found in
Kac' paper \cite{Kac:EGRG}.

Without giving details on the isomorphism, we note that the dual
of a bismash product can be viewed as a bismash product itself:
$$(k^G\# kF)^*\cong k^F\# kG.$$

\begin{Prop}
  Let $(F,G,\Hit,\Hitby)$ be a matched pair of finite groups.

  Then $\exp(k^G\# kF)=\exp(F\bowtie G)$.
\end{Prop}
\begin{proof}
  By a result of Beggs, Gould, and Majid \cite{BegGouMaj:FGFB},
  the Drinfeld double $D(k^G\# kF)$ is a Drinfeld twist of the
  Drinfeld double $D(k[G\bowtie F])$.
  By results of Etingof and Gelaki \cite{EtiGel:EFDHA} we have already cited,
  the exponent is invariant under such twists, and also invariant
  under taking the Drinfeld double.
\end{proof}

We shall now discuss the results of our computations of the
trivial power dimensions $\tpd_{m,n}(\QQ^G\# \QQ F)$ for some
matched pairs of groups. The computations were done with the help
of Maple; part of the code will be presented in \nmref{sec:code}.
In fact we have already presented results for special bismash
products, namely the Drinfeld doubles of some group algebras, in
the preceding section.

The next matched pair we consider is obtained from the symmetric
group $S_n$:
\begin{Lem}\nmlabel{Lemma}{Snfactorization}
  $S_n$ is factorizable into subgroups $F=S_{n-1}\subset S_n$
  (the copy of $S_{n-1}$ in $S_n$ that fixes $n$), and $G\cong C_n$,
  the cyclic group of order $n$ generated by the $n$-cycle
  $\tau=(1\,2\,3\dots n)$. Explicitly, for $\sigma\in S_n$ we have
  $\sigma=\sigma_1\sigma_2$
  with $\sigma_2=\tau^k\in G$, where $k=n-\sigma\inv(n)$, and
  $\sigma_1=\sigma\sigma_2\inv\in F$.
\end{Lem}
\begin{proof}
  It is easy to see that we have a factorizable group as stated,
  since the two subgroups have trivial intersection, and their
  orders multiply to give the order of $S_n$.
  If we want to factor $\sigma\in S_n$ into a product
  $\sigma_1\sigma_2$ with $\sigma_2=\tau^k\in G$ and $\sigma_1\in S_{n-1}$, we
  have to find an exponent $k$ such that $\sigma_1=\sigma\tau^{-k}$
  does not move $n$, that is, we have to make sure that
  $\tau^{-k}(n)=\sigma\inv(n)$, that is $n=\tau^k\sigma\inv(n)$, or
  simply $k=n-\sigma\inv(n)$.
\end{proof}
The matched pair arising from the factorizable group $S_n$ as in
the Lemma was studied in particular by Masuoka
\cite[Expl.1.3]{Mas:CSGHAE}, who shows that the Opext group for
this matched pair is trivial; in other words, the only middle term
of an extension \eqref{extpic} giving rise to this matched pair is
the bismash product.
\begin{Rem}
  The same formula as in the preceding Lemma also computes the
  factorization of an element $\sigma\in A_n$ for odd $n$ in the
  factorizable group $A_n\cong A_{n-1}\bowtie C_n$. The bismash
  product $k^{C_n}\# kA_{n-1}$ is naturally a Hopf subalgebra of
  $k^{C_n}\# kS_{n-1}$.
\end{Rem}

Tables~\ref{S4bismash} and \ref{dS4bismash} list the dimensions
$\tpd_{i,j}(H)$ for the bismash product $\QQ^{C_4}\# \QQ S_3$ and
its dual $\QQ^{S_3}\#\QQ C_4$, respectively. Note that we know the
exponent in each case to be that of $S_4$, namely $12$.
\begin{table}
$$ \begin {array}{|c|ccccccccccc|}
\hline
\tpd_{i,j}(\QQ^{C_4}\#\QQ S_3)&j=1&2&3&4&5&6&7&8&9&10&11\\
\hline i=1&1&1&1&1&1&1&1&1&1&1&1
\\
2&&12&1&8&5&12&3&12&5&8&1\\
%\hline
3& &&9&5&3&6&5&1&5&5&1\\
%\hline
4& & & &15&3&12&5&11&1&12&1
\\
%\hline
5& & & & &7&7&5&5&5&3&1\\
%\hline
6& && & & &18&7&12&6&12&1\\
%\hline
7& & & & & & &7&3&3&5&1\\
%\hline
8& & & & & & & &15&5&8&1\\
%\hline
9& & & & & & & & &9&1&1\\
%\hline
10& & & & & & & & & &12&1\\
%\hline
11& & & & & & & & & & &1
\\\hline
\end {array} $$
\caption{}\label{S4bismash}
\end{table}
\begin{table}
$$
\begin {array}{|c|ccccccccccc|}
\hline
  \tpd_{i,j}(\QQ^{S_3}\# \QQ C_4)&j=1&2&3&4&5&6&7&8&9&10&11
\\\hline
i=1&1&1&1&1&1&1&1&1&1&1&1\\
%\hline
2& &12&1&6&5&10&3&10&5&7&1\\
%\hline
3& & &9&5&3&9&5&1&5&5&1\\
%\hline
4& & &&15&3&10&5&11&1&10&1\\
%\hline
5& & & & &7&7&5&5&5&3&1\\
%\hline
6& & & & & &18&7&10&9&10&1\\
%\hline
7& & & & & & &7&3&3&5&1\\
%\hline
8& & & & & & & &15&5&6&1\\
%\hline
9& & & & & & & & &9&1&1\\
%\hline
10& & & & & & & & & &12&1\\
%\hline
11& & & & & & & & & & &1
\\\hline\end {array}
$$
\caption{}\label{dS4bismash}
\end{table} Using \nmref{recipe} we find:
\begin{Cor}
  The $24$-dimensional Hopf algebra $H=\QQ^{C_4}\# \QQ S_3$, of exponent
  $12$, contains elements of Hopf order $n$ for each $1\leq n\leq
  10$ except $n=7$ and $n=10$. In particular, it contains elements whose Hopf
  order is a prime that does not divide the dimension (or the
  exponent) of $H$.

  The dual Hopf algebra $H^*=\QQ^{S_3}\# \QQ C_4$ contains elements of
  Hopf order $n$ for each $1\leq n\leq 10$ except $n=7$ and $n=9$. In
  particular, the Hopf orders of elements that occur in $H$ and its
  dual are different.
\end{Cor}

\begin{Rem}\nmlabel{Remark}{indiv_orders}
While Tables~\ref{S4bismash} and \ref{dS4bismash} give us the full
information on which Hopf orders of elements are possible, they do
not tell us which individual elements have those orders. By a
separate calculation one can verify that
$$
\bas{(1\, 3)(2\, 4)}\#(2\, 3)-\bas{(1\, 4\, 3\, 2)}\# (1\,
2)+\bas{(1\, 3)(2\, 4)}\# (1\, 3\, 2)\in\QQ^{C_4}\#\QQ S_3=:H$$
has Hopf order five.

It is clear that all elements of the Hopf subalgebra
$\QQ^{C_4}\subset H$ have Hopf orders $1,2,$ or $4$. One can
verify that the elements of the subalgebra $\QQ S_3\subset H$ have
Hopf orders $1,2,6$ or $12$. While the last number may be
surprising, these possible orders are at least not prime to the
exponent of $H$, or of $\QQ S_3$. Even if elements of the form
$1\# a$ with $a\in S_3$ are, in this sense, well-behaved under the
Hopf power maps, their behavior does not seem to be easy to
understand and predict. For example, $1\# (1 3)$ has Hopf order
$2$, while $1\# (1 2)$ has Hopf order $12$. Note also that $3$ is
\emph{not} among the Hopf orders of elements in $\QQ S_3\subset
H$.

One can also verify that all elements of the standard basis $\bas
x\o a$ have Hopf orders $1,2,3,4,$ or $12$. Again, this does not
mean that the orders of such elements would be easy to understand.
For example, both $\bas{(1\,2\,3\,4)}\in\QQ^{C_4}$ and $(1\,2)\in
\QQ S_3$ have Hopf order $2$, but $\bas{(1\,2\,3\,4)}\#(1\,2)\in
H$ has Hopf order $3$. While $\bas{(1\,3)(2\,4)}\in\QQ^{C_4}$ has
Hopf order $4$, the element $\bas{(1\,3)(2\,4)}\#(1\,2)\in H$ has
Hopf order $2$. And $\bas{(1\,2\,3\,4)}\#(1\,2\,3)$ has Hopf order
$2$.
\end{Rem}

In a more indirect way, the example $\QQ^{C_4}\#\QQ S_3$ also
helps to answer a question raised in \nmref{sec:doubles} about the
behavior of Hopf powers under Drinfeld twists:

\begin{Expl}
  Consider the factorizable group $S_4$, and the bismash product
  $\QQ^{C_4}\# \QQ S_3$. From \nmref{grdbl} we know
  that $\tpd_n(D(kS_4))\neq 1$ if and only if $n$ is not prime to
  $4!$. But up to a Drinfeld twist, $D(kS_4)$ is isomorphic to
  $D(\QQ^{C_4}\# \QQ S_3)$ by \cite{BegGouMaj:FGFB}. On the other hand,
  we have seen that $\QQ^{C_4}\# \QQ S_3$, hence also its double,
  contains elements of Hopf order $5$. This shows that the
  property of having nontrivial elements with trivial $n$-th Hopf
  power is not invariant under twisting.
\end{Expl}

\begin{table}
$$% The table for the bismash product from the factorizable group a5
\begin {array}{|c|ccccccccccccccc|}
\cline{1-4}
\multicolumn{4}{|l|}{\tpd_{ij}(\QQ^{C_5}\# \QQ A_4)}\\
\cline{2-16}
&1&2&3&4&5&6&7&8&9&10&11&12&13&14&15\\\hline
1&\mathbf{   1 } &1&1&1&1&1&1&1&1&
1&1&1&1&1&1\\2&1&\mathbf{   32 }
&13&18&17&
18&9&20&1&18&11&26&11&12&21\\3&1&13&\mathbf{
37
 } &15&17&17&13&17&14&17&13&21&17&1&29\\4&1&18
&15&\mathbf{   34 } &7&20&13&22&13&26&13&22&11&16&19
\\5&1&17&17&7&\mathbf{   25 } &13&9&17&1&13&
13&17&9&1&25\\6&1&18&17&20&13&\mathbf{   36
 } &1&18&5&14&13&25&1&14&17\\7&1&9&13&13&9&1&
\mathbf{   21 }
&9&9&21&7&17&13&1&17\\8&1&20
&17&22&17&18&9&\mathbf{   32 } &1&22&9&18&9&12&17
\\9&1&1&14&13&1&5&9&1&\mathbf{   21 } &13&7&
13&9&1&17\\10&1&18&17&26&13&14&21&22&13&\mathbf{
  38 } &7&26&17&14&25\\11&1&11&13&13&13&
13&7&9&7&7&\mathbf{   19 }
&15&7&1&19\\12&1&
26&21&22&17&25&17&18&13&26&15&\mathbf{   44 } &13&10&29
\\13&1&11&17&11&9&1&13&9&9&17&7&13&\mathbf{   21
 } &1&21\\14&1&12&1&16&1&14&1&12&1&14&1&10&1&
\mathbf{   16 }
&1\\15&1&21&29&19&25&17&17&
17&17&25&19&29&21&1&\mathbf{   41 }
\\16&1&
12&1&16&1&14&1&12&1&14&1&10&1&16&1\\17&1&21&13&7&17&
9&9&13&1&9&11&17&11&1&21\\18&1&18&29&24&17&25&13&26&
13&26&13&30&17&10&29\\19&1&7&15&19&7&7&9&7&13&13&13&
13&11&1&19\\20&1&26&21&20&25&26&9&26&1&26&13&26&9&14
&25\\21&1&9&14&7&13&17&1&9&7&1&13&13&1&1&17
\\22&1&20&13&20&9&14&17&16&9&26&7&26&13&12&17
\\23&1&13&21&7&17&13&9&17&1&9&9&13&9&1&17
\\24&1&10&13&26&1&22&13&14&17&26&7&25&9&14&17
\\25&1&9&17&13&13&1&17&9&13&25&7&17&17&1&25
\\26&1&22&17&28&13&26&7&20&7&20&19&24&7&16&19
\\27&1&17&23&17&17&13&21&13&14&21&15&29&13&1&29
\\28&1&18&17&22&9&10&13&20&9&26&7&18&21&12&21
\\29&1&1&1&1&1&1&1&1&1&1&1&1&1&1&1\\\hline\end {array}
 $$
\caption{}\label{A5bismash}
\end{table}
\begin{table}
$$%the table for the dual of the bismash obtained from a5, i.e.
%\tpd_{ij}(k^{A_4}\# kC_5)
\begin {array}{|c|ccccccccccccccc|}
\cline{1-4}
\multicolumn{4}{|l|}{\tpd_{ij}(\QQ^{A_4}\# \QQ C_5)}\\
\cline{2-16}
&1&2&3&4&5&6&7&8&9&10&11&12&13&14&15\\
\hline1& \mathbf{ 1 } &1&1&1&1&1&1&1&1&
1&1&1&1&1&1\\2&1& \mathbf{ 32 } &13&22&17&
20&9&20&1&24&9&28&9&12&17\\3&1&13& \mathbf{ 37
 } &17&17&17&13&17&17&21&15&25&21&1&33\\4&1&22
&17& \mathbf{ 34 } &7&22&11&18&13&28&13&26&13&16&19
\\5&1&17&17&7& \mathbf{ 25 } &13&9&17&1&9&
13&17&9&1&21\\6&1&20&17&22&13& \mathbf{ 36
 } &1&20&9&16&13&28&1&16&21\\7&1&9&13&11&9&1&
 \mathbf{ 21 } &11&9&17&7&17&13&1&21\\8&1&
20&17&18&17&20&11& \mathbf{ 32 } &1&20&11&20&9&12&21
\\9&1&1&17&13&1&9&9&1& \mathbf{ 21 } &13&7&
17&9&1&21\\10&1&24&21&28&9&16&17&20&13& \mathbf{
38 } &7&28&21&16&21\\11&1&9&15&13&13&13&7&11&7
&7& \mathbf{ 19 } &13&7&1&19\\12&1&28&25&26
&17&28&17&20&17&28&13& \mathbf{ 44 } &13&12&33
\\13&1&9&21&13&9&1&13&9&9&21&7&13& \mathbf{ 21
 } &1&17\\14&1&12&1&16&1&16&1&12&1&16&1&12&1&
 \mathbf{ 16 } &1\\15&1&17&33&19&21&21&21&
21&21&21&19&33&17&1& \mathbf{ 41 } \\16&1&
12&1&16&1&16&1&12&1&16&1&12&1&16&1\\17&1&17&13&7&17&
13&9&13&1&9&9&17&9&1&17\\18&1&20&33&24&17&28&13&28&
17&28&15&32&17&12&33\\19&1&7&13&19&7&7&11&7&13&13&13
&15&9&1&19\\20&1&28&17&22&21&28&9&28&1&24&13&28&9&16
&21\\21&1&9&17&7&13&21&1&9&9&1&13&17&1&1&21
\\22&1&20&13&22&9&12&21&18&9&28&7&28&13&12&21
\\23&1&13&17&7&17&9&11&21&1&9&11&13&9&1&21
\\24&1&16&21&28&1&24&9&12&21&28&7&28&13&16&21
\\25&1&9&17&13&11&1&17&9&13&21&7&17&17&1&21
\\26&1&20&15&28&13&28&7&22&7&22&19&24&7&16&19
\\27&1&17&25&15&17&21&17&13&17&17&13&33&13&1&33
\\28&1&16&17&20&9&16&13&20&9&28&7&20&17&12&17
\\29&1&1&1&1&1&1&1&1&1&1&1&1&1&1&1\\\hline\end {array}
 $$
\caption{}\label{dA5bismash}
\end{table}
As we already saw for the double of $kA_4$, we have
$\tpd_{i,j}(H)\neq \tpd_{i,j}(H^*)$ in general. However, in the
last example as well as for the doubles studied in
\nmref{sec:doubles}, we see that the columns with prime index are
the same. Note, though, that this only holds for the printed parts
above the diagonal (otherwise the printed parts of the prime rows
would have to agree as well). Still, $\tpd_{pi}(H)=\tpd_{pi}(H^*)$
for $p$ prime and $i<p$ implies that those prime numbers $p$ such
that there is an element of Hopf order $p$ are the same for $H$
and $H^*$ (namely, $2,3$, and $5$). We will have to consider a
larger example to see that the last observation fails to be true
in general, and that we may have $\tpd_{pq}(H)\neq \tpd_{pq}(H^*)$
even for two primes $p,q$. Tables~\ref{A5bismash} and
\ref{dA5bismash} are for the bismash product $H=\QQ^{C_5}\# \QQ
A_4$ and its dual, associated to the factorizable group $A_5\cong
A_4\bowtie C_5$. The exponent of $A_5$ is $e=30$, and thus it is
preferable to use the symmetry \nmref{symmetry} to cut away
redundant parts of the table and save space. We have printed the
matrix $(\tpd_{i,j}(H))$ for column indices up to half the
exponent, i.e.\ we have printed the first 15 columns only. Since
$\tpd_{i,j}(H)=\tpd_{e-i,e-j}(H)$, the columns we cut off appear
upside down in reverse order in the printed part of the table,
that is, the $16$-th column is the $14$-th column upside down,
etc. We have boldfaced the diagonal elements to make them easier
to spot. Now we have the following recipe to decide if there is an
element of Hopf order $n$ in $H$. If $n\leq e/2$, check if the
$n$-th diagonal element is strictly larger than all the numbers in
the column above it; there exists an element of Hopf order $n$ if
and only if this is the case. If $n>e/2$, check if the $(e-n)$-th
diagonal element is strictly larger than all the numbers in the
column \emph{below} it; there exists an element of Hopf order $n$
if and only if this is the case. For example, $H=\QQ^{C_5}\# \QQ
A_4$ does not contain an element of Hopf order $17$; since
$17>15$, and $30-17=13$, this is checked by finding the $13$-th
diagonal element, which is $21$, in the thirteenth column two
places below the diagonal. Observe that the dual Hopf algebra
$H^*=\QQ^{A_4}\# \QQ C_5$ does contain an element of Hopf order
$17$, since $21$ does not appear below the diagonal in the $13$-th
column of the table for $H^*$. Thus we have an example where $H$
does not contain an element of a certain prime Hopf order, while
$H^*$ does. The same phenomenon occurs for the prime $13$, with
the roles of $H$ and $H^*$ reversed. Note also that
$\tpd_{3,13}(H)=17\neq 21=\tpd_{3,13}(H^*)$. Of course, one can
make a complete list of the orders of elements of $H$ and $H^*$.
In addition to $1$ and the exponent, which are always possible
orders, the candidates are the numbers in $\{2,\dots,e-2\}$. For
$H=\QQ^{C_5}\#\QQ A_4$, the orders $13$ and $21$ occur only in
$H$, while $17$ and $25$ occur only in $H^*$, and $14$, $16$,
$19$, $23$ occur in neither. The total number of possible orders
is the same, namely $22$, for both $H$ and $H^*$.

\begin{table}
 \scalebox{0.5}{$\begin {array}{|c|cccccccccccccccccccccccccccccc|} \hline
\tpd_{i,j}&1&2&3&4&5&6&
7&8&9&10&11&12&13&14&15&16&17&18&19&20&21&22&23&24&25&26&27&28&29&30
\\\hline 1& \mathbf{ 1 } &1&1&1&1&1&1&1&1&1&1&1&1&
1&1&1&1&1&1&1&1&1&1&1&1&1&1&1&1&1\\2&1& \mathbf{ 69 }
&17&45&25&47&17&51&9&53&11&61&11&47&25&40&29&51&15&59&17&54
&13&45&9&53&21&45&9&63\\3&1&17& \mathbf{ 57
 } &33&21&33&21&33&26&25&13&41&17&18&49&21&17&37&23&37&26&21&21&
33&17&29&43&37&5&45\\4&1&45&33& \mathbf{ 86
 } &15&53&21&65&21&55&13&70&11&43&37&60&15&53&27&70&21&51&7&70&
13&55&33&67&5&65\\5&1&25&21&15& \mathbf{ 37
 } &25&17&25&9&25&13&29&9&9&29&13&29&29&15&33&21&21&17&13&13&21&
21&17&13&37\\6&1&47&33&53&25& \mathbf{ 85
 } &13&59&13&51&13&73&1&49&33&51&21&66&19&63&25&51&13&61&1&63&29
&51&13&73\\7&1&17&21&21&17&13& \mathbf{ 33
 } &21&13&29&7&29&13&13&25&9&17&25&21&21&5&25&9&25&17&15&29&21&9
&33\\8&1&51&33&65&25&59&21& \mathbf{ 88 } &
17&55&9&72&9&55&33&56&21&63&19&78&25&52&17&68&9&59&29&70&9&67
\\9&1&9&26&21&9&13&13&17& \mathbf{ 41 } &13
&7&33&9&9&29&9&9&21&17&21&27&9&1&37&13&15&26&21&9&29
\\10&1&53&25&55&25&51&29&55&13& \mathbf{ 77
 } &7&64&17&44&33&51&21&59&25&61&9&64&9&56&25&55&29&59&5&73
\\11&1&11&13&13&13&13&7&9&7&7& \mathbf{ 19
 } &15&7&1&19&1&11&13&13&13&13&7&9&7&7&19&15&7&1&19
\\12&1&61&41&70&29&73&29&72&33&64&15& \mathbf{
104 } &13&53&49&64&29&69&25&78&33&64&13&85&17&63&49&70&13&81
\\13&1&11&17&11&9&1&13&9&9&17&7&13& \mathbf{ 21
 } &1&21&1&11&17&11&9&1&13&9&9&17&7&13&21&1&21
\\14&1&47&18&43&9&49&13&55&9&44&1&53&1& \mathbf{
69 } &21&35&9&49&13&59&10&43&1&53&1&50&17&47&9&53
\\15&1&25&49&37&29&33&25&33&29&33&19&49&21&21&
 \mathbf{ 61 } &21&25&37&27&45&29&25&17&37&25&31&49&41&5&49
\\16&1&40&21&60&13&51&9&56&9&51&1&64&1&35&21&
 \mathbf{ 76 } &13&43&13&58&21&47&1&60&1&43&17&63&5&53
\\17&1&29&17&15&29&21&17&21&9&21&11&29&11&9&25&13&
 \mathbf{ 33 } &25&15&25&17&21&13&13&9&19&21&19&13&33
\\18&1&51&37&53&29&66&25&63&21&59&13&69&17&49&37&43&
25& \mathbf{ 85 } &27&63&21&55&17&65&17&55&29&59&13&81
\\19&1&15&23&27&15&19&21&19&17&25&13&25&11&13&27&13&
15&27& \mathbf{ 31 } &19&11&21&7&25&13&21&21&23&9&31
\\20&1&59&37&70&33&63&21&78&21&61&13&78&9&59&45&58&
25&63&19& \mathbf{ 94 } &29&55&21&70&13&63&37&70&9&73
\\21&1&17&26&21&21&25&5&25&27&9&13&33&1&10&29&21&17&
21&11&29& \mathbf{ 41 } &9&9&25&1&21&26&17&9&29
\\22&1&54&21&51&21&51&25&52&9&64&7&64&13&43&25&47&21
&55&21&55&9& \mathbf{ 73 } &9&52&17&52&25&51&5&67
\\23&1&13&21&7&17&13&9&17&1&9&9&13&9&1&17&1&13&17&7&
21&9&9& \mathbf{ 21 } &1&9&13&13&9&1&21\\24
&1&45&33&70&13&61&25&68&37&56&7&85&9&53&37&60&13&65&25&70&25&52&1&
 \mathbf{ 96 } &13&55&37&66&13&73\\25&1&9&
17&13&13&1&17&9&13&25&7&17&17&1&25&1&9&17&13&13&1&17&9&13&
 \mathbf{ 25 } &7&17&17&1&25\\26&1&53&29&55&21&
63&15&59&15&55&19&63&7&50&31&43&19&55&21&63&21&52&13&55&7&
 \mathbf{ 77 } &27&53&9&65\\27&1&21&43&33&21&29&
29&29&26&29&15&49&13&17&49&17&21&29&21&37&26&25&13&37&17&27&
 \mathbf{ 57 } &33&5&45\\28&1&45&37&67&17&51&21&
70&21&59&7&70&21&47&41&63&19&59&23&70&17&51&9&66&17&53&33&
 \mathbf{ 86 } &5&63\\29&1&9&5&5&13&13&9&9&9&5&1
&13&1&9&5&5&13&13&9&9&9&5&1&13&1&9&5&5& \mathbf{ 13 } &13
\\30&1&63&45&65&37&73&33&67&29&73&19&81&21&53&49&53&
33&81&31&73&29&67&21&73&25&65&45&63&13& \mathbf{ 97 }
\\31&1&5&9&9&9&13&13&5&5&9&1&13&1&5&9&9&9&13&13&5&5&
9&1&13&1&5&9&9&9&13\\32&1&63&29&63&25&55&17&68&21&51
&11&78&11&52&37&56&29&51&15&78&29&51&13&62&9&57&33&65&9&63
\\33&1&21&49&29&25&33&17&33&34&25&13&41&17&18&41&21&
21&37&19&37&34&21&21&33&17&33&35&33&9&45\\34&1&49&31
&65&15&57&21&57&13&56&13&64&11&54&35&44&15&53&27&61&11&54&7&56&13&63&
33&57&5&65\\35&1&17&17&7&25&13&9&17&1&13&13&17&9&1&
25&1&17&17&7&25&13&9&17&1&13&13&17&9&1&25\\36&1&49&
37&70&25&73&13&66&25&53&13&85&1&46&37&70&21&65&19&70&37&53&13&82&1&56&
33&62&13&73\\37&1&9&13&13&9&1&21&9&9&21&7&17&13&1&17
&1&9&13&9&9&1&17&9&13&17&7&21&13&1&21\\38&1&54&25&49
&25&51&21&67&9&56&9&57&9&51&25&39&21&63&19&63&17&55&17&53&9&54&21&51&9
&67\\39&1&5&34&29&5&13&17&17&33&21&7&33&9&9&37&17&5&
21&21&21&19&17&1&37&13&11&34&29&5&29\\40&1&51&37&70&
25&59&29&70&21&73&7&78&17&52&45&70&21&59&25&82&21&63&9&70&25&61&37&78&
5&73\\41&1&23&17&21&25&25&15&21&15&19&19&27&7&13&23&
13&23&25&21&25&21&19&9&19&7&27&19&15&13&31\\42&1&59&
29&55&29&66&29&55&21&63&15&81&13&43&37&49&29&71&25&59&21&63&13&65&17&
53&37&51&13&81\\43&1&19&25&19&17&13&25&21&13&25&7&25
&21&13&29&9&19&29&23&21&5&21&9&21&17&15&21&29&9&33\\
44&1&44&17&62&9&45&13&65&21&43&1&62&1&53&21&54&9&49&13&70&17&39&1&70&1
&44&21&56&9&53\\45&1&29&41&33&33&33&21&33&37&33&19&
49&21&21&53&21&29&37&23&45&37&25&17&37&25&35&41&37&9&49
\\46&1&43&17&43&13&59&9&47&2&53&1&58&1&43&21&53&13&
43&13&52&9&51&1&46&1&54&18&52&5&53\\47&1&21&13&7&17&
9&9&13&1&9&11&17&11&1&21&1&21&13&7&17&9&9&13&1&9&11&17&11&1&21
\\48&1&53&49&70&29&73&25&78&33&61&13&90&17&58&49&62&
25&81&27&78&33&57&17&85&17&64&41&78&13&81\\49&1&7&15
&19&7&7&9&7&13&13&13&13&11&1&19&1&7&15&19&7&7&9&7&13&13&13&13&11&1&19
\\50&1&61&29&53&33&55&21&63&9&65&13&61&9&53&33&43&25
&63&19&73&21&56&21&53&13&56&25&51&9&73\\51&1&9&34&25
&17&25&9&17&19&9&13&33&1&2&37&21&13&21&15&21&33&9&9&25&1&13&34&21&5&29
\\52&1&51&29&67&21&59&25&70&17&63&7&78&13&47&33&65&
21&55&21&70&17&67&9&66&17&57&33&68&5&67\\53&1&21&25&
15&29&25&17&25&9&21&9&25&9&9&21&13&25&29&15&29&17&21&21&13&9&21&17&17&
13&33\\54&1&45&29&55&13&67&25&59&25&55&7&73&9&59&33&
45&13&66&25&59&13&51&1&73&13&57&33&55&13&73\\55&1&17
&25&21&21&13&29&21&17&33&7&29&17&13&33&9&17&29&25&25&5&25&9&25&25&15&
25&25&9&37\\56&1&49&29&74&21&55&15&67&25&53&19&70&7&
43&33&62&19&55&21&70&29&49&13&70&7&65&29&63&9&65\\57
&1&25&35&29&25&29&25&29&34&29&15&49&13&17&41&17&25&29&17&37&34&25&13&
37&17&31&49&29&9&45\\58&1&49&25&49&17&45&21&51&9&61&
7&53&21&43&29&44&19&59&23&51&5&54&9&49&17&49&21&63&5&63
\\59&1&1&1&1&1&1&1&1&1&1&1&1&1&1&1&1&1&1&1&1&1&1&1&1
&1&1&1&1&1&1\\\hline\end {array}
$
}
 \smallskip
 \caption{$\tpd_{i,j}=\tpd_{i,j}(\QQ^{C_5}\#\QQ S_4)$}\label{S5bism}
\end{table}
\begin{table}
 \scalebox{0.5}{$\begin {array}{|c|cccccccccccccccccccccccccccccc|}\hline
\tpd_{i,j}&1&2&3&4&5&6&
7&8&9&10&11&12&13&14&15&16&17&18&19&20&21&22&23&24&25&26&27&28&29&30
\\\hline1& \mathbf{ 1 } &1&1&1&1&1&1&1&1&1&1&1&1&
1&1&1&1&1&1&1&1&1&1&1&1&1&1&1&1&1\\2&1& \mathbf{ 69 }
&25&65&17&69&13&65&21&65&5&69&1&65&25&65&17&65&17&69&25&65&
5&65&1&69&25&61&13&69\\3&1&25& \mathbf{ 57
 } &31&21&37&21&37&27&29&15&45&21&21&53&21&17&41&21&37&27&25&17&
41&17&29&45&35&5&45\\4&1&65&31& \mathbf{ 86
 } &13&69&21&72&31&77&7&80&11&69&31&76&13&69&23&76&21&73&1&86&11
&69&29&78&13&77\\5&1&17&21&13& \mathbf{ 37
 } &21&17&29&9&13&13&29&9&13&25&13&29&21&15&31&21&13&17&13&11&21
&21&21&13&21\\6&1&69&37&69&21& \mathbf{ 85
 } &13&73&29&69&9&77&1&69&37&69&21&81&17&77&37&65&9&77&1&77&37&
65&13&85\\7&1&13&21&21&17&13& \mathbf{ 33
 } &21&13&21&7&29&13&13&29&13&17&17&23&21&5&21&11&21&17&13&25&25
&9&22\\8&1&65&37&72&29&73&21& \mathbf{ 88
 } &21&65&9&80&9&65&37&72&25&73&17&88&29&61&17&72&9&73&33&80&13&
73\\9&1&21&27&31&9&29&13&21& \mathbf{ 41 }
&29&7&37&9&21&31&21&9&33&17&21&29&29&1&41&13&21&27&27&9&37
\\10&1&65&29&77&13&69&21&65&29& \mathbf{ 77
 } &5&73&9&69&29&69&13&69&21&69&21&73&1&77&9&69&29&69&13&77
\\11&1&5&15&7&13&9&7&9&7&5& \mathbf{ 19 } &
13&7&1&19&1&9&13&13&11&13&5&11&7&7&9&13&7&1&15\\12&1
&69&45&80&29&77&29&80&37&73&13& \mathbf{ 104 } &13&65&53&72&
29&77&27&86&37&73&13&88&17&69&53&80&13&85\\13&1&1&21
&11&9&1&13&9&9&9&7&13& \mathbf{ 21 } &1&17&1&9&5&9&9&1&5&9&13
&17&1&13&15&1&9\\14&1&65&21&69&13&69&13&65&21&69&1& 65&1& \mathbf{
69 } &21&69&13&65&13&69&21&65&1&69&1&69&21&65&
13&69\\15&1&25&53&31&25&37&29&37&31&29&19&53&17&21&  \mathbf{ 61 }
&21&21&41&27&39&31&29&21&41&21&29&53&35&5&49
\\16&1&65&21&76&13&69&13&72&21&69&1&72&1&69&21&
 \mathbf{ 76 } &13&65&13&76&21&65&1&76&1&69&21&72&13&69
\\17&1&17&17&13&29&21&17&25&9&13&9&29&9&13&21&13&
 \mathbf{ 33 } &17&15&31&17&13&13&13&9&21&25&21&13&21
\\18&1&65&41&69&21&81&17&73&33&69&13&77&5&65&41&65&
17& \mathbf{ 85 } &21&73&37&65&9&77&5&73&37&69&13&85
\\19&1&17&21&23&15&17&23&17&17&21&13&27&9&13&27&13&
15&21& \mathbf{ 31 } &19&11&21&7&25&13&17&23&19&9&27
\\20&1&69&37&76&31&77&21&88&21&69&11&86&9&69&39&76&
31&73&19& \mathbf{ 94 } &31&65&17&76&9&77&39&80&13&77
\\21&1&25&27&21&21&37&5&29&29&21&13&37&1&21&31&21&17
&37&11&31& \mathbf{ 41 } &21&9&29&1&29&27&21&9&37
\\22&1&65&25&73&13&65&21&61&29&73&5&73&5&65&29&65&13
&65&21&65&21& \mathbf{ 73 } &1&73&9&65&29&65&13&73
\\23&1&5&17&1&17&9&11&17&1&1&11&13&9&1&21&1&13&9&7&
17&9&1& \mathbf{ 21 } &1&9&9&13&9&1&10\\24&
1&65&41&86&13&77&21&72&41&77&7&88&13&69&41&76&13&77&25&76&29&73&1&
 \mathbf{ 96 } &13&69&37&78&13&85\\25&1&1&
17&11&11&1&17&9&13&9&7&17&17&1&21&1&9&5&13&9&1&9&9&13& \mathbf{ 25
 } &1&17&15&1&9\\26&1&69&29&69&21&77&13&73&21&
69&9&69&1&69&29&69&21&73&17&77&29&65&9&69&1& \mathbf{ 77 } &
29&65&13&77\\27&1&25&45&29&21&37&25&33&27&29&13&53&
13&21&53&21&25&37&23&39&27&29&13&37&17&29& \mathbf{ 57 } &33&
5&45\\28&1&61&35&78&21&65&25&80&27&69&7&80&15&65&35&
72&21&69&19&80&21&65&9&78&15&65&33& \mathbf{ 86 } &13&69
\\29&1&13&5&13&13&13&9&13&9&13&1&13&1&13&5&13&13&13&
9&13&9&13&1&13&1&13&5&13& \mathbf{ 13 } &13
\\30&1&69&45&77&21&85&22&73&37&77&15&85&9&69&49&69&
21&85&27&77&37&73&10&85&9&77&45&69&13& \mathbf{ 97 }
\\31&1&13&9&13&9&13&13&13&5&13&1&13&1&13&9&13&9&13&
13&13&5&13&1&13&1&13&9&13&9&13\\32&1&69&33&72&27&69&
21&80&21&65&7&86&9&65&35&72&27&65&19&86&27&65&13&72&9&69&35&76&13&69
\\33&1&25&47&31&25&37&17&37&37&29&15&45&21&21&43&21&
21&41&17&37&37&25&17&41&17&29&35&35&9&45\\34&1&65&29
&77&13&69&21&65&29&77&5&73&9&69&29&69&13&69&21&69&21&73&1&77&9&69&29&
69&13&77\\35&1&5&17&1&25&9&9&17&1&1&13&17&9&1&21&1&
17&9&7&19&13&1&17&1&11&9&17&9&1&9\\36&1&69&37&76&25&
85&13&80&29&69&13&88&1&69&41&76&25&81&19&86&41&65&9&84&1&77&41&72&13&
85\\37&1&1&13&9&9&1&21&9&9&9&7&17&13&1&21&1&9&5&11&9
&1&9&11&9&17&1&17&13&1&10\\38&1&65&29&65&21&73&13&73
&21&65&9&65&1&65&29&65&17&73&17&73&29&61&9&65&1&73&25&65&13&73
\\39&1&21&37&31&5&29&17&21&31&29&7&37&9&21&41&21&5&
33&21&21&19&29&1&41&13&21&37&27&5&37\\40&1&65&39&86&
21&69&29&80&31&77&7&88&19&69&39&76&21&69&23&84&21&73&9&86&19&69&37&86&
13&77\\41&1&17&19&19&25&21&15&21&15&17&19&25&7&13&23
&13&21&25&21&23&21&17&11&19&7&21&17&19&13&27\\42&1&
69&37&73&17&77&21&65&37&73&9&85&5&65&41&65&17&77&25&69&33&73&5&81&9&69
&41&65&13&85\\43&1&13&29&23&17&13&25&21&13&21&7&25&
21&13&25&13&17&17&21&21&5&17&9&25&17&13&21&27&9&21\\
44&1&65&21&76&13&69&13&72&21&69&1&72&1&69&21&76&13&65&13&76&21&65&1&76
&1&69&21&72&13&69\\45&1&25&43&31&29&37&25&37&41&29&
19&53&17&21&51&21&25&41&23&39&41&29&21&41&21&29&43&35&9&49
\\46&1&65&21&69&13&69&13&65&21&69&1&65&1&69&21&69&13
&65&13&69&21&65&1&69&1&69&21&65&13&69\\47&1&5&13&1&
17&9&9&13&1&1&9&17&9&1&17&1&21&5&7&19&9&1&13&1&9&9&21&9&1&9
\\48&1&65&53&78&29&81&25&88&37&69&15&92&17&65&53&72&
25&85&25&88&37&65&17&88&17&73&45&86&13&85\\49&1&5&13
&11&7&5&11&5&13&9&13&15&9&1&19&1&7&9&19&7&7&9&7&13&13&5&15&7&1&15
\\50&1&69&29&69&21&77&13&73&21&69&9&69&1&69&29&69&21
&73&17&77&29&65&9&69&1&77&29&65&13&77\\51&1&25&37&21
&17&37&9&29&19&21&13&37&1&21&41&21&13&37&15&31&31&21&9&29&1&29&37&21&5
&37\\52&1&65&33&80&21&65&29&76&29&73&5&88&13&65&37&
72&21&65&21&80&21&73&9&80&17&65&37&80&13&73\\53&1&17
&21&13&29&21&19&29&9&13&11&25&9&13&25&13&25&21&15&29&17&13&21&13&9&21&
17&21&13&22\\54&1&65&37&77&13&77&21&65&37&77&5&81&9&
69&37&69&13&77&21&69&29&73&1&85&9&69&37&69&13&85\\55
&1&13&25&23&19&13&29&21&17&21&7&29&17&13&29&13&17&17&25&21&5&21&9&25&
25&13&25&27&9&21\\56&1&69&29&76&23&77&13&80&21&69&11
&78&1&69&31&76&23&73&19&86&31&65&9&76&1&77&31&72&13&77
\\57&1&25&35&29&25&37&21&33&37&29&13&53&13&21&43&21&
29&37&19&39&37&29&13&37&17&29&47&33&9&45\\58&1&61&25
&69&13&65&17&65&25&69&5&65&5&65&25&65&13&69&17&65&21&65&1&69&5&65&25&
69&13&69\\59&1&1&1&1&1&1&1&1&1&1&1&1&1&1&1&1&1&1&1&1
&1&1&1&1&1&1&1&1&1&1\\\hline\end {array}$
}
 \smallskip
 \caption{$\tpd_{i,j}=\tpd_{i,j}(\QQ^{S_4}\#\QQ C_5)$}\label{S5bismdu}
\end{table}
We have already seen one example (the double of $S_4$) in
\nmref{sec:doubles} where a Hopf algebra and its dual admit
different total numbers of Hopf orders. Another example is the
bismash product of the factorizable group $S_5$ in
Tables~\ref{S5bism} and \ref{S5bismdu}. While $45$ of the numbers
in $\{2,\dots,58\}$ are Hopf orders of elements of $\QQ^{C_5}\#
\QQ S_4$, only $31$ are Hopf orders of elements of $\QQ^{S_4}\#
\QQ C_5$. While in this last example the percentage of possible
Hopf orders that actually occur is rather small compared to the
other examples in this section (let alone the doubles of small
groups in Tables~\ref{DS3}, \ref{DA4}, and \ref{dDA4}, where
\emph{all} conceivable orders occurred), it is still much larger
than in group algebras.

Although the picture might change when larger dimensions are
considered, it seems quite unlikely that further computations will
lead us to examples where a prime $p$ divides the exponent of a
bismash product Hopf algebra $H$, but $p$ is not the Hopf order of
an element of $H$. Note that if $H=k^G\# kF$, then $p$, a divisor
of the order of $F\bowtie G$, has to divide the order of either
$F$ or $G$; then $p$ occurs as the Hopf order of an element of one
of the Hopf algebras $k^F$ and $k^G$. These are Hopf subalgebras
in $H^*$ and $H$, respectively, so that $p$ is at least the Hopf
order of an element of $H$ \emph{or} $H^*$. Since the trivial
power dimensions are invariant under taking the dual, this also
implies that nontrivial elements with trivial $p$-th Hopf power
exist both in $H$ and $H^*$.

According to Tables~\ref{A5bismash}, \ref{dA5bismash},
\ref{S5bism} and \ref{S5bismdu} the relevant primes $2,3,5$ occur
as Hopf orders both in $H$ and $H^*$ for $H=\QQ^{C_5}\# \QQ S_4$
as well as $H=\QQ^{C_5}\# \QQ A_4$. As we have already seen in
\nmref{indiv_orders}, however, it is far from obvious how to find
elements of a given Hopf order $n$ in $k^G\# kF$ if one knows
elements of order $n$ in $F$. The possible Hopf orders we found
for elements of the the subalgebra $\QQ S_4=\QQ\neut\#\QQ
S_4\subset \QQ^{C_5}\# \QQ_4$ are $1, 2, 4, 12,$ and $ 30$. In
particular, the subalgebra $\QQ A_4\subset \QQ^{C_5}\# \QQ A_4$
also contains no elements of Hopf order $3$, since we can
naturally view $\QQ^{C_5}\# \QQ A_4\subset \QQ^{C_5}\# \QQ S_4$ as
a Hopf subalgebra. The element
$$\bas{(1\, 2\, 3\, 4\, 5)}\#(1\, 2\, 4)\in \QQ^{C_5}\# \QQ A_4\subset
\QQ^{C_5}\# \QQ S_4$$ has Hopf order $3$.

\section{Maple code}\nmlabel{Section}{sec:code}

We will now present the Maple programming used to compute the
numerical results that we have discussed in the preceding
sections.

For the moment, we assume that we have already implemented in
Maple a matched pair $(F,G,\Hit,\Hitby)$ of finite groups,
including explicit bijections $E_M\colon M\to\{1,\dots,|M|\}$ for
$M=F,G$. We will talk later about how to actually provide the
following items:
\def\itt#1:{\item[{\tt #1}]}
\begin{description}
  \item[{\tt NumfromF}] A maple procedure implementing $E_F$.
  \item[{\tt FfromNum}] A maple procedure implementing $E_F\inv$.
  \item[{\tt ordF}] The order of $F$.
  \itt multF: A maple procedure that computes the product of
  two elements of $F$.
  \itt oneF: The neutral element of $F$.
  \itt invF: A maple procedure that computes the inverse of an
  element in $F$.
  \itt NumfromG: The same as {\tt NumfromF} for the group $G$. We
  also have {\tt GfromNum}, {\tt ordG}, {\tt multG}, {\tt invG},
  {\tt oneG}.
  \itt hit: The procedure call {\tt hit}$(x,a)$ should compute
  $x\Hit a$.
  \itt hitby: The procedure call {\tt hitby}$(x,a)$ should compute
  $x\Hitby a$.
  \itt expL: The exponent of the group $L=F\bowtie G$.
\end{description}

From these data we proceed to compute {\tt ordL:=ordF*ordG} and a
procedure to enumerate the elements of the standard basis of
$H=k^G\#kF$. The element $\bas x\# a$ will be represented in Maple
as a two-element list {\tt [x,a]}. The following procedures
enumerate the elements in the standard basis in inverse
lexicographic order:
\begin{tt}
\begin{verbatim}
> NumfromB:=proc(h)
>    option remember;
>    (NumfromF(h[2])-1)*ordG+NumfromG(h[1]);
> end proc;
>
> BfromNum:=proc(i)
>    option remember;
>    [ GfromNum(irem(i-1,ordG)+1),
>      FfromNum(iquo(i-1,ordG)+1) ];
> end proc;
\end{verbatim}
\end{tt}

Based on these enumeration procedures, we will represent a general
element of $H$ as a column vector of length {\tt ordL}. To handle
vectors and matrices, we use Maple's {\tt LinearAlgebra} package.
The result of multiplying two basis elements will either be zero
or a basis element. We represent the result as a general element,
that is, a column vector:
\begin{tt}
\begin{verbatim}
> multBB:=proc(h,j) option remember;
>    x:=h[1]; a:=h[2]; y:=j[1]; b:=j[2];
>    if hitby(x,a)=y
>       then UnitVector(NumfromB([x,multF(a,b)]),ordL);
>       else ZeroVector(ordL);
>    end if;
> end proc;
\end{verbatim}
\end{tt}
Besides multiplying two basis elements, we also need to be able to
multiply a general element of $H$ and an element of the standard
basis. This is of course easily reduced to the multiplication of
basis elements.
\begin{tt}
\begin{verbatim}
> multHB:=proc(V,h);
>    R:=ZeroVector(ordL); # a register to sum into
>    for i to ordL do
>        if not V[i]=0
>           then B:=V[i]*multBB(BfromNum(i),h); # the i-th summand
>        end if;
>        R:=R+B; # is added to the register
>    end do;
>    R; # in the end, the register contains the result.
> end proc;
\end{verbatim}
\end{tt}

Note that $\Delta(\bas x\# a)$ is a sum of simple tensors that
happen to be tensor products of two elements of the standard
basis. Thus, we can represent the result of comultiplication on a
basis element as a list of two-element lists of basis elements.
This is computed by the following procedure:
\begin{tt}
\begin{verbatim}
> comult:=proc(h);
>    x:=h[1]; a:=h[2]; # so h=[x,a]
>    [seq([[multG(x,invG(GfromNum(i))),
>           hit(GfromNum(i),a)],
>          [GfromNum(i),a]],
>     i=1..ordG)];
> end proc;
\end{verbatim}
\end{tt}
The map $[n]\colon H\to H$ is naturally represented in Maple by
its matrix with respect to the standard basis, in the convention
that the matrix for $f\colon H\to H$ with respect to a basis $h_i$
is the matrix $m_{ij}$ with $f(h_j)=\sum m_{ij}h_i$, that is, the
$i$-th column vector of the matrix is the coordinate vector of the
image under $f$ of the $i$-th basis vector. The map $[1]$ is the
identity, represented by the unit matrix. To compute $[n]$
recursively, we use the formula $h\pow {n+1}=h\sw 1\pow nh\sw 2$.
We assume we are given the matrix $A$ representing the $n$-th
power map, and we wish to compute the matrix representing the
$(n+1)$-st power map. Recall that $h\sw 1\o h\sw 2$, for a basis
element $h$, is represented as a list of two-element lists. For
each element in that list, which represents a simple tensor $f\o
g$, we compute the image of $f$ under the $n$-th power map, simply
by looking up the relevant column in the matrix $A$. Then we
multiply the element represented by that column and the basis
element $g$, using the procedure {\tt multHB} given above. This we
do for each simple tensor $f\o g$ in the list which is the output
of the comultiplication procedure applied to $h$, and sum up the
results. This gives one column of the desired matrix for the next
power map.
\begin{tt}
\begin{verbatim}
> NextPowerMatrix:=proc(A);
>    # The parameter A is assumed to be the matrix representing
>    # the n-th Hopf power endomorphism of H. The result of
>    # NextPowerMatrix should be the matrix representing the
>    # n+1-st Hopf power endomorphism.
>    R:=Matrix(ordL,ordL); # a register to compute
>                          # the resulting matrix in
>    for i to ordL do   # run through the basis of H
>        T:=comult(BfromNum(i));  # comultiply basis element
>        LengthofT:=ordG;         # to get a list of this length
>        C:=ZeroVector(ordL);     # a register to compute
>                                 # the i-th column of the result
>        for j to LengthofT do    # go through the summands
>            D:=A[1..ordL,NumfromB(T[j][1])];
>                   # look up the column in the previous matrix
>                   # corresponding to the first tensor factor in
>                   # the summand under consideration. So D represents
>                   # the previous Hopf power of that tensor factor.
>            C:=C+multHB(D,T[j][2]);  # multiply this with the second
>                                     # tensor factor, and add to
>                                     # the register.
>        end do;
>        R[1..ordL,i]:=C;  # store this in the i-th column.
>    end do;
>    R;
> end proc;
\end{verbatim}
\end{tt}

To compute all the matrices, say $A_1,\dots,A_e$ for the $i$-th
power maps, we invoke
\begin{tt}
\begin{verbatim}
> A[1]:=IdentityMatrix(ordL);
> for l from 2 to expL do
>     A[l]:=NextPowerMatrix(A[l-1])
> end do;
\end{verbatim}
\end{tt}
The matrix for the map $\eta\epsilon\colon H\to H$, which maps
$\bas x\# a$ to $\sum_y\bas y\#\neut$ if $x=\neut$, and to zero
otherwise, is computed by the following procedure:
\begin{tt}
\begin{verbatim}
> etaepsilon:=Matrix(ordL,ordL,
>   (i,j) -> if BfromNum(j)[1]=oneG and BfromNum(i)[2]=oneF
>               then 1
>               else 0
>            end if);
\end{verbatim}
\end{tt}
We can then compute bases for the spaces $\TPS_{n}$ by
\begin{tt}
\begin{verbatim}
> for i to expL-1 do
>     T[i]:=NullSpace(A[i]-etaepsilon)
> end do;
\end{verbatim}
\end{tt}
and the dimensions $\tpd_{i,j}(H)$ as the number of elements in a
basis for the intersection of two such spaces by
\begin{tt}
\begin{verbatim}
> t[i,j]:=nops(IntersectionBasis([T[i],T[j]]));
\end{verbatim}
\end{tt}
If we want to know the Hopf order of a specific element of $H$
(represented by a vector of length the order of $L$), we can
compute it, given the matrices for the power maps as above, by
\begin{tt}
\begin{verbatim}
> HopfOrder:=proc(h)
>    # computes the Hopf order of an element
>    i:=1; # the least possible Hopf order is 1
>    while not Equal(A[i].h,etaepsilon.h) do
>       i:=i+1 # while the i-th power is not trivial, add one
>    end do;
>    i; # so this is the least i for which the i-th power is trivial
> end proc;
\end{verbatim}
\end{tt}
In \nmref{sec:bismash} we gave an example of an element of Hopf
order $5$ in $H=k^{C_4}\# kS_3$. This was found more or less by
trial and error using the procedure {\tt HopfOrder}. The elements
of the basis {\tt T[5]} computed by Maple for the fifth trivial
power space were found to have Hopf orders $2$ and $3$. The sum of
a basis element of Hopf order $2$ and one of Hopf order $3$ had
Hopf order $5$.

 To deal with $H^*$, it is sufficient to observe that the
matrices for the $n$-th power map and $\eta\epsilon$ for $H^*$ can
be found by transposing the matrices computed for $H$.

Now we should go about providing the numbers and procedures
necessary to deal with the groups $F,G$, which we had assumed to
be given above.

First, we implement the factorizable group $S_n\cong
S_{n-1}\bowtie C_n$. We wrote down the factorization of $\sigma\in
S_n$ as a product $\sigma=\sigma_1\sigma_2$ for unique
$\sigma_1\in S_{n-1}$ and $\sigma_2=(1\,2\dots n)^k\in C_n$ in
\nmref{Snfactorization}. The three procedures {\tt Snfactor2exp},
{\tt Snfactor2}, and {\tt Snfactor1} below compute $k,\sigma_2$,
and $\sigma_1$, respectively, given $\sigma$ and the rank $n$ of
the symmetric group as their two arguments. They are based on
Maple's {\tt group} package. As auxiliary procedures we provide
{\tt applyinvperm}, which applies the inverse of its first
argument, a permutation, to its second, a number, and {\tt
ncyclepower} which provides the powers of the standard $n$-cycle,
given $n$ as its first, and the exponent as its second argument.
We provide and use a procedure {\tt multperms} for multiplying
permutations in the order used in this paper (acting on the left
of elements), reversing the convention in Maple's {\tt group}
package.
\begin{tt}
\begin{verbatim}
> Snfactor2exp:=proc(sigma,n)
>    option remember;
>    n-applyinvperm(sigma,n);
> end proc;
>
> Snfactor2:=proc(sigma,n)
>    option remember;
>    ncyclepower(n,Snfactor2exp(sigma,n));
> end proc;
>
> Snfactor1:=proc(sigma,n)
>    option remember;
>    multperms(sigma,invperm(Snfactor2(sigma,n)));
> end proc;
>
> applyinvperm:=proc(sigma,i);
>    # applies the inverse of permutation sigma to element i
>    j:=i;  # the result will be j
>    for cyc in sigma do # for each cycle (a list)
>        if member(i,cyc,'k') # look if it contains i,
>                       # and remember at which position
>           then if k=1 # if in the first position
>                then j:=cyc[nops(cyc)] # the inverse cycle
>                          # maps i to the last element
>                else j:=cyc[k-1] # otherwise to the predecessor
>           end if
>        end if
>    end do;
>    j;
> end proc;
>
> ncyclepower:=proc(n,e);
>    convert([seq(i+ e mod n +1,i=0..n-1)],'disjcyc');
> end proc;
>
> multperms:=proc(x,y);
>    mulperms(y,x);
> end proc;
\end{verbatim}
\end{tt}

Now we fix a Maple variable {\tt N} for the rank of a symmetric
group (in practice, $N=4,5$ were practicable and interesting), and
set up a matched pair $(F,G,\Hit,\Hitby)$ with $F\cong S_{N-1}$
and $G\cong C_N$. Recall that the operations of the two groups on
each other are obtained from a factorization in $S_N$ of the
product $xa$ of $x\in G$ and $a\in F$.
\begin{tt}
\begin{verbatim}
> hit:=proc(x,a)
>    Snfactor1(multperms(x,a),N);
> end proc;
>
> hitby:=proc(x,a)
>    Snfactor2(multperms(x,a),N);
> end proc;
>
\end{verbatim}
\end{tt}
We also need all the other procedures for handling $F$ and $G$
that we assumed to exist above:
\begin{tt}
\begin{verbatim}
> FList:=elements(permgroup(N-1,{ncycle(N-1),[[1,2]]}));
>   # To set up a bijection between elements of F and
>   # numbers, we get a list of all elements of F.
>   # Conversion then consists in looking up elements
>   # in the list.
>
> ordF:=(N-1)!
> NumfromF:=proc(a)
>    member(a,FList,'i');i;
> end proc;
>
> FfromNum:=proc(i)
>    FList[i];
> end proc;
>
> multF:=proc(a,b)
>    # Multiplication in F is that of permutations.
>    multperms(a,b);
> end proc;
>
> invF:=proc(a)
>    # Inverse in F is that for permutations.
>    invperm(a);
> end proc;
>
> oneF:=[];
>
> ordG:=N;
>
> NumfromG:=proc(x)
>    # The number we assign to a power of the N-cycle is
>    # the exponent, also the distance it moves 1.
>    applyperm(x,1);
> end proc;
>
> GfromNum:=proc(n)
>    ncyclepower(N,n-1);
> end proc;
>
> multG:=proc(s,t)
>    mulperms(t,s);
> end proc;
>
> invG:=proc(s)
>    invperm(s);
> end proc;
>
> oneG:=[];
>
\end{verbatim}
\end{tt}
Since we are only dealing with a few ranks $N$, it is much easier
to simply enter the exponent {\tt expL} by hand than to devise a
Maple procedure to compute it, so {\tt expL:=12} for $N=4$, and
{\tt expL:=60} for $N=5$.

The same {\tt hit} and {\tt hitby} procedures that stem from the
factorization of $S_N$ can equally well be used for the
factorization of $A_N$, for odd $N$ (in practice, $N=5$) into
subgroups isomorphic to $A_{N-1}$ and $C_N$, respectively. The
group $G$ is the same as above, and the procedures for $F$ only
change where we set up a list of elements of $F$ to enumerate the
elements. For $N=5$ we use
\begin{tt}
\begin{verbatim}
> Flist:=elements(permgroup(4,{[[1, 2, 3]], [[2, 3, 4]]}));
> expL:=30;
\end{verbatim}
\end{tt}
to produce a list of all the elements of $A_4$.

Next, let us discuss how to do our computations for the Drinfeld
double of a group $G$: We set up numbers and procedures {\tt
NumfromG, GfromNum, ordG, multG, invG} as before, and use the same
procedures again for the group $F=G$. The exponent of the group
$L\cong G\times G$ is the same as the exponent of $G$. The action
$\Hit$ should be trivial, and the action $\Hitby$ should be the
adjoint action of $G$ on itself from the right:
\begin{tt}
\begin{verbatim}
> hit:=proc(x,a) a;end proc;
>
> hitby:=proc(x,a);
>    option remember;
>    multG(invG(a),multG(x,a));
> end proc;
\end{verbatim}
\end{tt}
The computations for the Hopf algebra $k^G\o kG$, from which the
double $D(kG)$ is obtained by a dual Drinfeld twist, can be done
by defining both actions to be trivial. Computations for a group
algebra or its dual can be done by defining one of the two groups
$F,G$ as well as both actions to be trivial. Of course this is not
particularly efficient, and was really done by different Maple
procedures. Also, the computations for $k^G\o kG$ can be done more
efficiently by computing the power matrices for $kG$ first, and
then taking their Kronecker products with their respective
transposes, to obtain the power matrices for $k^G\o kG$. We have
omitted listing these additional procedures to keep the paper to a
reasonable length.
\bibliographystyle{acm}
\bibliography{eigene,andere}
\end{document}